\documentclass{mscs}
\usepackage{enumerate,hyperref,tikz,amssymb,xspace,amsmath}
\usepackage[all]{xy}
\usetikzlibrary{matrix,arrows}
\DeclareMathOperator{\Spec}{Spec}
\DeclareMathOperator{\fin}{Fin}
\DeclareMathOperator{\compfin}{KFin}
\DeclareMathOperator{\Span}{Span}
\DeclareMathOperator{\Proj}{Proj}
\DeclareMathOperator{\colim}{colim}
\newcommand{\ie}{\textit{i.e.}\xspace}

\newcommand{\op}{\ensuremath{^{\text{op}}}}
\newcommand{\A}{\ensuremath{\mathcal{A}}}
\newcommand{\AF}{\ensuremath{\mathrm{AF}}}
\newcommand{\Id}{\ensuremath{\mathrm{Idl}}}
\newcommand{\B}{\ensuremath{\mathcal{B}}}
\newcommand{\C}{\ensuremath{\mathbb{C}}}
\newcommand{\CC}{\ensuremath{\mathcal{C}}}
\newcommand{\D}{\ensuremath{\mathcal{D}}}
\newcommand{\F}{\ensuremath{\mathcal{F}}}
\newcommand{\G}{\ensuremath{\mathcal{G}}}
\newcommand{\HH}{\ensuremath{\mathcal{H}}}
\newcommand{\K}{\ensuremath{\mathcal{K}}}
\newcommand{\N}{\ensuremath{\mathbb{N}}}
\newcommand{\U}{\ensuremath{\mathcal{U}}}
\newcommand{\V}{\ensuremath{\mathcal{V}}}
\newcommand{\up}{\ensuremath{\mathop{\uparrow}}}
\newcommand{\down}{\ensuremath{\mathop{\downarrow}}}
\newcommand{\ddown}{\ensuremath{\mathop{\rotatebox[origin=c]{270}{$\twoheadrightarrow$}}}}
\newcommand{\uup}{\ensuremath{\mathop{\rotatebox[origin=c]{90}{$\twoheadrightarrow$}}}}
\newcommand{\OMP}{\ensuremath{\mathbf{OMP}}}

\newcommand{\DCPO}{\ensuremath{\mathbf{DCPO}}}
\newcommand{\CStar}{\ensuremath{\mathbf{CStar}}}

\numberwithin{equation}{section}
\newtheorem{theorem}[equation]{Theorem}
\newtheorem{lemma}[equation]{Lemma}
\newtheorem{proposition}[equation]{Proposition}
\newtheorem{corollary}[equation]{Corollary}
\newtheorem{definition}[equation]{Definition}
\newtheorem{example}[equation]{Example}

\begin{document}

\title{Domains of commutative C*-subalgebras}

\author[C. Heunen and B. Lindenhovius]
  {C\ls H\ls R\ls I\ls S\ns H\ls E\ls U\ls N\ls E\ls N$^1$%
  \ns and\ns B\ls E\ls R\ls T\ns L\ls I\ls N\ls D\ls E\ls N\ls H\ls O\ls V\ls I\ls U\ls S$^2$\\
  $^1$ School of Informatics, University of Edinburgh, \addressbreak 10 Crichton Street, Edinburgh EH8 9AB, UK.
  \addressbreak $^2$ Department of Computer Science, Tulane University, \addressbreak 6823 St.\ Charles Ave, New Orleans LA 70118, US.}

\maketitle

\begin{abstract}
  \noindent
  A C*-algebra is determined to a great extent by the partial order of its commutative C*-subalgebras.
  We study order-theoretic properties of this dcpo.
  Many properties coincide: the dcpo is, equivalently, algebraic, continuous, meet-continuous, atomistic, quasi-algebraic, or quasi-continuous, if and only if the C*-algebra is scattered. 
  For C*-algebras with enough projections, these properties are equivalent to finite-dimensionality. 
  Approximately finite-dimensional elements of the dcpo correspond to Boolean subalgebras of the projections of the C*-algebra.
  Scattered C*-algebras are finite-dimensional if and only if their dcpo is Lawson-scattered. 
  General C*-algebras are finite-dimensional if and only if their dcpo is order-scattered. 
\end{abstract}

\section*{Introduction}

One can study a C*-algebra $A$ through the partially ordered set $\CC(A)$ of its commutative C*-subalgebras.
In general this poset does not determine the C*-algebra completely, which follows from Connes' example of a C*-algebra $A$ that is not isomorphic to its opposite C*-algebra $A^\text{op}$, \ie\ the algebra with the same underlying vector space as $A$ but  multiplication $(a,b)\mapsto ba$~\cite{connes:opposite}, for which clearly $\CC(A)\simeq\CC(A^\text{op})$. However, $\CC(A)$ does determine the structure of $A$ to a great extent:
\begin{itemize}
  \item $\CC(A)$ determines $A$ up to quasi-Jordan isomorphism~\cite{hamhalter:ordered};
  \item $\CC(A)$ determines $A$ up to a Jordan isomorphism if $A$ is an AW*-algebra\footnote{This statement and the next one were originally proved under the condition that the algebras do not have without type I$_2$ summands. This condition can be removed~\cite[Corollary 9.2.9]{lindenhovius:thesis}.}~\cite{hamhalter:dye}; 
  \item $\V(A)$, a variant of $\CC(A)$ defined below in Definition~\ref{def:AofA}, determines $A$ up to Jordan isomorphism if $A$ is a W*-algebra~\cite{doeringharding:jordan}; 
  \item $\CC(A)$ determines $A$ up to $*$-isomorphism if $A$ is a type I AW*-algebra~\cite[Corollary 8.6.24]{lindenhovius:thesis};
  \item $\CC(A)$, together with extra structure making it a so-called \emph{active lattice}, determines $A$ up to $*$-isomorphism if $A$ is an AW*-algebra~\cite{heunenreyes:activelattice}.   
\end{itemize}
Thus $\CC(A)$ can be used as a substitute for the C*-algebra itself~\cite{heunen:manyfaces}.

The intuition is clearest in the case of quantum theory. There, the C*-algebra models all observations one can possibly perform on a quantum system. However, not all observations may be performed simultaneously, but only those that live together in a commutative C*-subalgebra $C$.
There is an inherent notion of approximation: if $C \subseteq D$, then $D$ contains more observations, and hence provides more information. 

This sort of informational approximation is more commonly studied in \emph{domain theory}~\cite{abramskyjung:domaintheory,gierzetal:domains}. 
In domain theory, elements of a poset are sometimes interpreted as incomplete objects which are missing information. The commutative C*-subalgebras of a C*-algebra do not contain information as such in the direct sense. Instead, the incompleteness resides in the way the information in the quantum system is observed. In this sense the poset $\CC(A)$ is more akin to the poset of the domain of definition of partial functions, ordered by extension. Classically this is not a very interesting poset, but in the quantum setting it is.
As we are speaking of a continuous amount of observables, but in practice only have access to a discrete number of them, we are most interested in partial orders where every element can be approximated by empirically accessible ones. Domain theory has a variety of notions modelling this intuition. 

This article studies how these domain-theoretic properties of $\CC(A)$ relate to operator-algebraic properties of the C*-algebra $A$. We show that they all coincide in our setting, as the following are equivalent:
\begin{itemize}
  \item the C*-algebra $A$ is scattered (as defined in Section~\ref{sec:cstaralgebras});
  \item the partial order $\CC(A)$ is algebraic (Section~\ref{sec:algebraicity});
  \item the partial order $\CC(A)$ is continuous (Section~\ref{sec:continuity});
  \item the partial order $\CC(A)$ is meet-continuous (Section~\ref{sec:meetcontinuity});
  \item the right adjoint of $\CC(f)$ is Scott-continuous for each injective $*$-homomorphism $f \colon B \to A$ (Section~\ref{sec:meetcontinuity});
  \item the partial order $\CC(A)$ is atomistic (Section~\ref{sec:atomicity});
  \item the partial order $\CC(A)$ is quasi-algebraic (Section~\ref{sec:quasi});
  \item the partial order $\CC(A)$ is quasi-continuous (Section~\ref{sec:quasi}).
\end{itemize}
This makes precise exactly `how much approximate finite-dimensionality' on the analytical side is required for these desirable notions of approximation on the domain-theoretic side. It is satisfying that these notions robustly coincide with the established algebraic notion of scatteredness, which is intimately related to approximate finite-dimensionality. 

We also study finite-dimensionality of $A$ in terms of the partial order $\CC(A)$:
\begin{itemize}
  \item a C*-algebra $A$ is finite-dimensional if and only if $\CC(A)$ is Lawson-scattered (Section~\ref{sec:scattered});
  \item a C*-algebra $A$ is finite-dimensional if and only if $\CC(A)$ is order-scattered (Section~\ref{sec:orderscattered});
\end{itemize}
Finally, we study the links between domain theory and projections, which form an important part of traditional C*-algebra theory:
\begin{itemize}
  \item the partial order $\CC_\AF(A)$ of commutative approximately finite-dimensional C*-sub\-algebras of a C*-algebra $A$ is isomorphic to the domain of Boolean subalgebras of $\Proj(A)$, the projections in $A$, which allows us to reconstruct $\Proj(A)$ from $\CC(A)$ (Section~\ref{sec:Booleansubalgebras}).
  \item an AW*-algebra $A$ is finite-dimen\-sional if and only if $\CC(A)$ is continuous, if and only if $\CC(A)$ is algebraic (Section~\ref{sec:projections});
  \item the functors $\CC$ and $\CC_\AF$ do not preserve directed colimits of C*-algebras, whereas the functor $\B$ that assigns to each orthomodular poset its poset of Boolean subalgebras does preserve directed colimits of orthomodular posets (Section~\ref{sec:colimits}).
\end{itemize}

\section{C*-algebras}\label{sec:cstaralgebras}

For the benefit of readers with a background in domain theory, we briefly recall what we need from the classical theory of C*-algebras~\cite{kadisonringrose:oa1,weggeolsen:ktheory,takesaki:oa1,conway:functionalanalysis}.

\begin{definition}
  A \emph{norm} on a complex vector space $V$ is a function $\|-\| \colon V \to [0,\infty)$ satisfying
  \begin{itemize}
    \item $\|v\|=0$ if and only if $v=0$;
    \item $\|\lambda v\|=|\lambda|\|v\|$ for $\lambda \in \C$;
    \item $\|v+w\|\leq\|v\|+\|w\|$.
  \end{itemize}
  A \emph{Banach space} is a normed vector space that is complete in the metric $d(v,w)=\|v-w\|$.
\end{definition}

\begin{definition}
  An \emph{inner product} on a complex vector space $V$ is a map $\langle - \mid -\rangle \colon V \times V \to \C$ that:
  \begin{itemize}
    \item is linear in the second variable;
    \item is conjugate symmetric: $\langle v \mid w \rangle =\overline{\langle w \mid v \rangle}$; 
    \item satisfies $\langle v \mid v\rangle\geq 0$ with equality only when $v=0$.
  \end{itemize}
  An inner product space $V$ is a \emph{Hilbert space} when the norm $\|v\|=\sqrt{\langle v \mid v \rangle}$ makes it a Banach space.
  For example, $\C^n$ with its usual inner product $\langle v\mid w\rangle=\sum_{i=1}^n\bar v_iw_i$ for $v=(v_1,\ldots,v_n)$ and $w=(w_1,\ldots,w_n)$ is a Hilbert space.
\end{definition}

\begin{definition}
  A complex vector space $A$ is an \emph{algebra} when it carries a bilinear associative multiplication $A \times A \to A$. It is called \emph{unital} when it has a unit $1 \in A$ satisfying $1a=a=a1$.
  It is \emph{commutative} when $ab=ba$ for all $a,b \in A$.
  A $*$-algebra is an algebra $A$ with an \emph{involution}, i.e., a map $A \to A$, $a\mapsto a^*$, satisfying for each $a,b\in A$ and each $\lambda,\mu\in\C$:
  \begin{itemize}
  	\item $(\lambda a+\mu b)^*=\bar\lambda a^*+\bar\mu b^*$;
    \item $(a^*)^*=a$;
    \item $(ab)^* = b^* a^*$.
  \end{itemize}
  A \emph{C*-algebra} is a $*$-algebra $A$ that is simultaneously a Banach space with for each $a,b\in A$:
  \begin{itemize}
    \item $\|ab\| \leq \|a\| \|b\|$;
    \item $\|a^*a\| = \|a\|^2$.
  \end{itemize}
\end{definition}
The last identity is called the \emph{C*-identity}. One can show that the norm of a C*-algebra is completely determined by the algebraic structure of the algebra, for which the C*-identity is crucial. As a consequence every C*-algebra has a unique norm (see for instance \cite[Corollary C.28]{landsman:foundations}).
We emphasise that we will always assume our C*-algebras to be unital. This is due to the fact that Theorem~\ref{thm:hamhalter} below, which is fundamental for most results in this contribution, is only known to hold in the unital case. Similarly $*$-homomorphisms, which we will define below, are always assumed to preserve the unit.
An element $a$ of a C*-algebra is \emph{self-adjoint} when $a=a^*$, and a \emph{projection} when $a^*=a=a^2$. We write $\Proj(A)$ for the set of all projections in $A$.

\begin{example}\label{ex:BH}
  As mentioned, the set of all $n$-by-$n$ complex matrices is a C*-algebra, with the involution given by the conjugate transpose. More generally, the space $B(H)$ of all \emph{bounded operators} on a Hilbert space $H$, i.e., all continuous linear maps $a \colon H \to H$, is a C*-algebra as follows.
  Addition and scalar multiplication are defined by $a+b \colon v \mapsto a(v)+b(v)$ and $\lambda a\colon v\mapsto \lambda a(v)$, multiplication is composition by $ab \colon v \to a(b(v))$, and $1$ is the identity map $v \mapsto v$.
  The involution is defined by taking the \emph{adjoint}: given a bounded operator $a \colon H \to H$, we let $a^*$ be the unique bounded operator satisfying $\langle v \mid a(w) \rangle = \langle a^*(v) \mid w \rangle$ for each $v,w\in H$.
  The norm is given by $\|a\|=\sup\{\|a(v)\| \mid v \in H, \|v\|=1\}$. Notice that this C*-algebra is noncommutative (unless $H$ is one-dimensional or zero-dimensional). Moreover, $B(H)$ equals the algebra of all $n$-by-$n$ complex matrices if we choose $H=\C^n$.
\end{example}

The previous example is in fact prototypical, as the following theorem shows, for which we first need to introduce the appropriate morphisms of C*-algebras.
A linear map $f \colon A \to B$ between C*-algebras is a (unital) \emph{$*$-homomorphism} when $f(ab)=f(a)f(b)$, $f(a^*)=f(a)^*$. It is \emph{unital} when $f(1)=1$; in this article all $*$-homomorphisms are assumed to be unital. If $f$ is bijective we call it a \emph{$*$-isomorphism}, and write $A \simeq B$. We denote the category of C*-algebras and $*$-homomorphisms by $\CStar$, and note that the isomorphisms in this category are precisely the $*$-isomorphisms. Every $*$-homomorphism is automatically continuous, and is even an isometry when it is injective \cite[Theorem 4.1.8]{kadisonringrose:oa1}. Since any $*$-homomorphism $f\colon A\to B$ is linear, continuity of $f$ implies that its operator norm $\|f\|$ satisfies $\|f(a)\|\leq\|f\|\|a\|$ for each $a\in A$ \cite[Theorems 1.5.5 \& 1.5.6]{kadisonringrose:oa1}. A C*-algebra $B$ is a \emph{C*-subalgebra} of a C*-algebra $A$ when $B \subseteq A$, and the inclusion $B \to A$ is a $*$-homomorphism. As a consequence $B$ must contain the identity element of $A$. Moreover, since the inclusion must be an isometry, it follows that every C*-subalgebra of $A$ is a closed subset of $A$, and conversely, every norm-closed *-subalgebra of $A$ is a C*-subalgebra of $A$. Clearly, the inverse image $f^{-1}[B]$ of a $*$-homomorphism $f:A\to B$ is a C*-subalgebra of $A$. A less trivial fact is that the image $f[A]$ of $f$ is a C*-subalgebra of $B$ \cite[Theorem 4.1.9]{kadisonringrose:oa1}. 
\begin{theorem}[Gelfand--Naimark]\label{thm:gns}\cite[Theorem 4.5.6 \& Remark 4.5.7]{kadisonringrose:oa1}.
  Any C*-algebra is $*$-iso\-mor\-phic to a C*-subalgebra of $B(H)$ for a Hilbert space $H$.
\end{theorem}
The above C*-algebra is noncommutative unless $H$ is one-dimensional or zero-dimensional. Here is an example of a commutative one.

\begin{example}
  The vector space $\C^n$ is a commutative C*-algebra under pointwise operations in the max norm $\|(x_1,\ldots,x_n)\| = \max \{|x_1|,\ldots,|x_n|\}$. It sits inside the algebra $B(\C^n)$ of $n$-by-$n$ matrices as the subalgebra of diagonal ones, illustrating Theorem~\ref{thm:gns}. 
  \end{example}
  The infinite version of the previous example is as follows.
\begin{example}  Write $C(X)$ for the set of all continuous functions $f \colon X \to \C$ on a compact Hausdorff space $X$.
  It becomes a commutative C*-algebra as follows: addition and scalar multiplication are pointwise, \textit{i.e.}, $f+g \colon x \mapsto f(x)+g(x)$, multiplication is pointwise $fg \colon x \mapsto f(x)g(x)$, the unit is the function $x \mapsto 1$, the involution is given by $f^* \colon x \mapsto \overline{f(x)}$, and the norm is $\|f\| = \sup_{x \in X} |f(x)|$.
\end{example}

The above example is prototypical for commutative C*-algebras.

\begin{theorem}[Gelfand duality]\label{thm:gelfandduality}\cite[Theorem 4.4.3]{kadisonringrose:oa1}.
  Any commutative (unital) C*-al\-ge\-bra is $*$-isomorphic to $C(X)$ for a compact Hausdorff space $X$
  called its \emph{spectrum}.
\end{theorem}
\begin{proof}[Proof (sketch)]
  If $A$ is a commutative unital C*-algebra, define its Gelfand spectrum to be the set $X$ of all nonzero\footnote{Instead of demanding $\varphi$ be nonzero, we can require $\varphi(1)=1$. In fact, one can prove that $\varphi$ is actually a unital  $*$-homomorphism $A\to\C$.} linear maps $\varphi:A\to\mathbb{C}$ such that $\varphi(ab)=\varphi(a)\varphi(b)$. This is a subset of the unit ball of the dual $A^*$ of $A$, \ie\ the space of bounded functionals $A\to\C$. The dual $A^*$ becomes a Hausdorff space when equipped with the weak*-topology, which is generated by a subbasis consisting of sets of the form $\{\psi\in A^* \mid |\varphi(a)-\psi(a)|<\varepsilon\}$ where $\varphi\in A^*$, $a\in A$ and $\varepsilon>0$. Compactness of $X$ follows from the Banach-Alaoglu Theorem. For each $a\in A$, there is a continuous function $\hat a\colon X\to\mathbb{C}$ defined by $\hat a(\varphi)=\varphi(a)$. The *-isomorphism $A\to C(X)$, called the \emph{Gelfand transform}, is defined by $a\mapsto \hat a$. 
\end{proof}

The previous theorem extends to a categorical duality, which was first shown explicitly in \cite{negrepontis:duality}, see also \cite[Theorem C.23]{landsman:foundations}. We only need the following proposition, which guarantees that studying the poset of C*-subalgebras of a commutative C*-algebra reduces to studying compact Hausdorff quotients of its Gelfand spectrum.

\begin{proposition}\label{prop:dualitysubalgebrasandquotients}
  Let $A$ be a commutative C*-algebra with spectrum $X$. If $X\to Y$ is a continuous surjection onto a compact Hausdorff space $Y$, then $Y$ is homeomorphic to the spectrum of a C*-subalgebra of $A$. 
  Conversely, if a C*-subalgebra of $A$ has spectrum $Y$, there is a continuous surjection $X\to Y$. 
\end{proposition}
\begin{proof}
  If $q \colon X \to Y$ is a continuous surjection, then $B = \{ f \circ q \mid f\in C(Y) \}$ is a C*-subalgebra of $A$. Conversely, if $B$ is a C*-subalgebra of $A$, define a equivalence relation $\sim_B$ on $X$ by setting $x \sim_B y$ if and only if $b(x)=b(y)$ for each $b \in B$. The quotient $Y=X / {\sim_B}$ is a compact Hausdorff space and it follows that $C(Y)$ is *-isomorphic to $B$. For details, see~\cite[Proposition 5.1.3]{weaver:mathematicalquantization}.
\end{proof}

It is easy to see that the intersection of any collection of C*-subalgebras of a C*-algebra $A$ is again a C*-subalgebra of $A$. Hence if $S$ is a subset of a C*-algebra $A$, there is a smallest C*-subalgebra $C^*(S)$ of $A$ that contains $A$, which we call the \emph{C*-subalgebra of $A$ generated by $S$}. If $S=\{a_1,\ldots,a_n\}$ is finite, we write $C^*(a_1,\ldots,a_n)$ instead of $C^*(\{a_1,\ldots,a_n\})$. If $S$ consists of mutually commuting elements and is closed under $a\mapsto a^*$, then it can be embedded into a commutative *-subalgebra of $A$
, whose closure is a commutative C*-subalgebra, and so $C^*(S)$ is commutative. C*-subalgebras generated by a subset behave well under $*$-homomorphisms.
 	\begin{lemma}\label{lem:imageofgeneratedsubalgebra}
 		Let $f \colon A\to B$ be a $*$-homomorphism between C*-algebras $A$ and $B$, and let $S\subseteq A$ be a subset. Then $f[C^*(S)]=C^*(f[S])$.
 	\end{lemma}
 	\begin{proof}
 		Clearly $f[S]\subseteq f[C^*(S)]$ and so $C^*(f[S])\subseteq f[C^*(S)]$.
 		For the other inclusion, note $$S\subseteq f^{-1}\big[f[S]\big]\subseteq f^{-1}[C^*(f[S])].$$ Since the inverse image of a C*-subalgebra under a $*$-homomorphism is clearly a C*-subalgebra, it follows that $f^{-1}[C^*(f[S])] \subseteq A$ is a C*-subalgebra. Hence $C^*(S)\subseteq f^{-1}[C^*(f[S])]$, and $f[C^*(S)]\subseteq C^*(f [S])$.
 	\end{proof}

\subsection*{Commutative C*-subalgebras}

We now come to our main object of study, namely commutative C*-subalgebras. 
When generally describing (quantum) systems C*-algebraically, observables become self-adjoint elements $a=a^* \in A$. 
For each self-adjoint element $a$, there is a unique injective $*$-homomorphisms $C(\sigma(a)) \to A$ sending function $x\mapsto x$ to $a$~\cite[Theorem 4.4.5]{kadisonringrose:oa1} linking observables to commutative C*-subalgebras. 
Here $\sigma(a)$ is the \emph{spectrum} of $a$, i.e., the compact Hausdorff space
\[
  \sigma(a) = \{ \lambda \in \C \mid a - \lambda 1 \text{ is not invertible}\}.
\]
The following definition captures the main structure of approximation on the algebraic side.

\begin{definition}
  For a C*-algebra $A$, define
  \[  
    \CC(A) = \{ C \subseteq A \mid C \text{ is a commutative C*-subalgebra}\},
  \]
  partially ordered by inclusion: $C \leq D$ when $C \subseteq D$.
\end{definition}

Let us consider some elementary domain-theoretic properties of $\CC(A)$ now. For detailed information about domain theory, we refer to~\cite{abramskyjung:domaintheory,gierzetal:domains}.
Let $\CC$ be a partially ordered set. We think of its elements as partial computations or observations, and the partial order $C \leq D$ as ``$D$ provides more information about the eventual outcome than $C$''. 
With this interpretation, it is harmless to consider \emph{downsets}, or \emph{principal ideals}, instead of $C \in \CC$:
\[
  \down C = \{D\in\CC \mid D\leq C\}.
\]
Dually, it is also of interest to consider \emph{upsets}, or \emph{principal filters}, consisting of all possible expansions of the information contained in $C \in \CC$:
\[
  \up C = \{D\in\CC \mid D\geq C\}.
\]
This extends to subsets $\D\subseteq\CC$ as:
\[
  \down \D = \bigcup_{D\in\D}\down D, 
  \qquad
  \up\D = \bigcup_{D\in\D}\up D.
\]
If $\D$ has a least upper bound in $\CC$, it is denoted by $\bigvee\D$. 
Furthermore, $\D$ is called \emph{directed} if for each $D_1,D_2\in\D$ there is a $D_3\in\D$ such that $D_1,D_2\leq D_3$.
This can be interpreted as saying that the partial computations or observations in $\D$ can always be compatibly continued without leaving $\D$. If we want to emphasize that the set $\D$ over which we take the supremum if directed, we write $\bigvee^{\up}\D$ instead of $\bigvee\D$.
Similarly, we write $\bigwedge \D$ for a greatest lower bound, when it exists.
For two-element sets $\D$ we just write the \emph{meet} $\bigwedge \{D_1, D_2\}$ as $D_1 \wedge D_2$.

\begin{definition}
  A partially ordered set $\CC$ is \emph{directed-complete partially ordered} (dcpo) if each directed subset of $\CC$ has a least upper bound.
\end{definition}

\begin{proposition}\label{prop:CAisdcpo}
  If $A$ is any C*-algebra, then $\CC(A)$ is a dcpo, where the supremum $\bigvee\D$ of a directed set $\D\subseteq\CC(A)$ is given by $\overline{\bigcup\D}$.
\end{proposition}
\begin{proof}
   Let $\D\subseteq\CC(A)$ be a directed subset. Let $S=\bigcup\D$. We show that $S$ is a commutative *-algebra. Let $x,y\in S$ and $\lambda,\mu\in\C$, there are $D_1,D_2\in\D$ such that $x\in D_1$ and $y\in D_2$. Since $\D$ is directed, there is some $D_3\in\D$ such that $D_1,D_2\subseteq D_3$. Hence $x,y\in D_3$, whence $\lambda x+\mu y,x^*,xy\in D_3$, and since $D_3$ is commutative, it follows that $xy=yx$. Since $D_3\subseteq S$, it follows that $S$ is a commutative *-subalgebra of $A$. Then $\overline{S}$ is a commutative C*-subalgebra of $A$, which is clearly the least upper bound of $\D$. 
  See also~\cite{spitters:partialmeasurement}.
\end{proof}

In order to show that the assignment $A \mapsto \CC(A)$ is also functorial, we first have to introduce the appropriate notion of morphisms of dcpos.
A function $f:P\to Q$ between partially ordered sets $P$ and $Q$ is \emph{monotone} when $p \leq q$ in $P$ implies $f(p) \leq f(q)$ in $Q$; it is an \emph{order embedding} if it is monotone and $f(p)\leq f(q)$ implies $p\leq q$, and it is is an \emph{order isomorphism} is a monotone bijection with a monotone inverse, or equivalently, if it is a surjective order embedding. If $P$ and $Q$ are dcpos, then we say that a monotone map $f:P\to Q$ is \emph{Scott-continuous} function it it preserves the suprema of directed subsets. The category of dcpos with Scott continuous maps is denoted by $\DCPO$.

The next proposition shows that $\CC$ is a functor $\CStar\to\DCPO$.
\begin{proposition}\label{prop:CfisScottcontinuous}
Let $f:A\to B$ be a $*$-homomorphism between C*-algebras $A$ and $B$. Then the map $\CC(f):\CC(A)\to\CC(B)$, $C\mapsto f[C]$ is Scott continuous. In particular, if $f$ is injective, then $\CC(f)$ is an order embedding. 
\end{proposition}
 \begin{proof}
Let $f:A\to B$ be a $*$-homomorphism, and let $C\subseteq A$ be a commutative C*-subalgebra. Since the image of a $*$-homomorphism is a C*-subalgebra of the codomain, it follows that $f[C]$ is a C*-subalgebra of $B$. Since $C$ is commutative, and $f$ preserves all algebraic operations, it follows that $f[C]$ is commutative. Hence the assignment $C\mapsto f[C]$ is a well-defined map $\CC(A)\to\CC(B)$, which clearly preserves inclusions, hence it is monotone. Let $\D\subseteq\CC(A)$ be directed. 
  Since $\CC(f)$ is monotone, $\CC(f)[\D]$ is directed. 
  Now 
  \[
  f \left[\bigvee\D \right]
  = f \left[\overline{\bigcup\D}\right]
  \subseteq \overline{f\left[\bigcup\D\right]} 
  \subseteq \overline{f\left[\overline{\bigcup\D}\right]} = \overline{f\left[\bigvee\D\right]}
  =  f\left[\bigvee\D\right],
  \]
  where the first inclusion holds because $*$-homomorphisms are continuous, and where the last equality holds because C*-subalgebras are closed. 
  Thus 
  \[
f\left[\bigvee\D\right]=\overline{f\left[\bigcup\D\right]}  =\overline{\bigcup_{D\in\D}f[D]} = \bigvee_{D\in\D}f[D] ,
  \]
  hence $\CC(f)$ is Scott continuous. Finally, let $f$ be injective and assume that $\CC(f)(C)\subseteq\CC(f)(D)$, i.e., $f[C]\subseteq f[D]$. Let $x\in C$. Then $f(x)\in f[D]$, so there is some $y\in D$ such that  $f(x)=f(y)$. By injectivity of $f$ it follows that $x=y$, hence $x\in D$. We conclude that $\CC(f)$ is an order embedding. 
\end{proof}

The dcpo $\CC(A)$ is of interest because it determines the C*-algebra $A$ itself to a great extent, as mentioned in the Introduction. The following theorem, which generalizes an earlier result in the setting of compact Hausdorff quotients of compact Hausdorff spaces~\cite[Theorem 11]{mendivil:functionalgebras}, illustrates this.

\begin{theorem}\cite[Theorem 2.4]{hamhalter:ordered}\label{thm:hamhalter}.
  Let $A,B$ be commutative C*-algebras. Given any order isomorphism $\psi \colon \CC(A) \to \CC(B)$, there exists a $*$-iso\-mor\-phism $f \colon A \to B$ such that $\CC(f)=\psi$ that is unique unless $A$ is two-dimensional.
\end{theorem}

If $A$ is two-dimensional, the $*$-isomorphism $A \to B$ is not unique because $\mathbb C^2$ has two automorphisms, the identity and the flip map $(x,y)\mapsto(y,x)$, that both induce the identity automorphism on $\CC(\mathbb C^2)$.

It already follows that for arbitrary C*-algebras $A$, the partial order on $\CC(A)$ determines the C*-algebra structure of each individual element of $\CC(A)$. Indeed, if $C\in\CC(A)$, then $\down C$ is order isomorphic to $\CC(C)$, and since $C$ is a commutative C*-algebra, it follows that the partially ordered set $\down C$ determines the C*-algebra structure of $C$. 

\subsection*{Approximate finite-dimensionality}

In practice, within finite time one can only measure or compute up to finite precision, and hence can only work with (sub)systems described by finite-dimensional C*-subalgebras.
Therefore one might think that the natural extension is for the finite-dimensional C*-subalgebras to be dense in the whole C*-algebra. C*-algebras that can be described in this way are called \emph{approximately finite-dimensional}. In the separable case these can be classified in several ways, for instance by means of Bratteli diagrams~\cite{bratteli:afalgebras} or by K-theory~\cite{weggeolsen:ktheory}. 

\begin{definition}\label{def:af}
 We call a C*-algebra $A$:
  \begin{itemize}
  	\item  \emph{approximately finite-dimensional}, or an \emph{AF-algebra}, if there is a directed set $\D$ of finite-dimen\-sional C*-subalgebras of $A$ whose union is dense in $A$ with respect to the norm topology.
  	\item  \emph{locally approximately finite-dimensional}, or a \emph{locally AF-algebra}, when for each $\varepsilon>0$ and each $a_1,\ldots,a_n\in A$ there exist a finite-dimensional C*-subalgebra $B\subseteq A$ and $b_1,\ldots,b_n\in B$ such that $\|a_i-b_i\|<\varepsilon$ for any $i=1,\ldots,n$. 
  \end{itemize}
\end{definition}
Before we give some examples, we need to develop some theory, but the reader could already fast-forward to Example \ref{ex:cantor} for an example of a commutative AF-algebra, and to Example \ref{ex:compact operators} for an example of a non-commutative locally AF-algebra. Let us point out that we do not, as most authors do, restrict approximately finite-dimensional C*-algebras to be separable (in which case the directed set $\D$ can be replaced by a chain of finite-dimensional C*-algebras), since all our results also hold in the non-separable case. 
It is easy to see that any AF-algebra is locally approximately finite-dimensional: if $a_1,\ldots, a_n$ and $\varepsilon>0$, then there are $b_1,\ldots,b_n\in\bigcup\D$ such that $\|a_i-b_i\|<\varepsilon$. Each $a_i$ is contained in some $D_i\in\D$, and since $\D$ is directed, there is some $B\in\D$ containing $D_1,\ldots,D_n$, hence also $\{b_1,\ldots,b_n\}$. In case $A$ is separable, the converse holds \cite[Theorem 2.2]{bratteli:afalgebras}, but in general the class of locally AF-algebras is strictly larger than the class of AF-algebras~\cite{farahkatsura:dixmier}. However, Proposition \ref{prop:commutativeAFistotallydisconnected} below shows that for \emph{commutative} C*-algebras both notions coincide. Moreover, it turns out that a commutative C*-algebra is approximately finite-dimensional if and only if its spectrum is \emph{totally disconnected}: that is, when its connected components are exactly the singletons. Separability gives the additional requirement that the spectrum be second-countable~\cite[Proposition 3.1]{bratteli:structurespacesofafalgebras}. The proposition is well known, and we provide a proof for convenience. We first need a lemma.

\begin{lemma}\label{lem:finiteprojectionsgeneratefinitedimsubalgebra}
Let $X$ be compact Hausdorff, and let $p_1,\ldots,p_n$ be projections in $C(X)$. Then $C^*(p_1,\ldots,p_n)$ is a finite-dimensional subalgebra of $C(X)$, and is spanned by all finite products of elements in the set $\{p_1,\ldots,p_n,1\}$.
\end{lemma}	
\begin{proof}
	Let $S$ be the collection of all finite products of $p_1,\ldots,p_n$ and $1$. Since $C(X)$ is commutative, the $p_i$ mutually commute, and moreover, since they are idempotent, it follows that $S$ must be finite. As a consequence, $S$ is closed under the multiplication. Since all projections are self-adjoint and commute, their products should be self-adjoint, too, hence $S$ is closed under the involution $a\mapsto a^*$. Now let $V$ be the span of $S$. Then $V$ is finite-dimensional, and since $S$ is closed under the multiplication and the involution, it follows that $V$ is a *-subalgebra of $C(X)$, which is finite-dimensional, hence topologically closed. Thus $V$ is a C*-subalgebra of $C(X)$ containing the projections $p_1,\ldots,p_n$, hence it must contain $C^*(p_1,\ldots,p_n)$. Since the latter must contain $S$, hence also $V$, we find that $V=C^*(p_1,\ldots,p_n)$. We conclude that $C^*(p_1,\ldots,p_n)$ must be finite-dimensional.
\end{proof}

\begin{proposition}\label{prop:commutativeAFistotallydisconnected}
  The following are equivalent for a compact Hausdorff space $X$:
  \begin{itemize}
  \item[(1)] $C(X)$ is approximately finite-dimensional;
  \item[(2)] $C(X)$ is locally approximately finite-dimensional;
  \item[(3)] $X$ is totally disconnected.
  \item[(4)] $C(X)$ is generated by its projections;
  \end{itemize}
\end{proposition} 
\begin{proof}
	We already showed that a C*-algebra is locally approximately finite-dimensional if it is approximately finite-dimensional, which yields (1)$\implies$(2).  For (2)$\implies$(3), assume that $C(X)$ is locally approximately finite-dimensional, and let $x,y\in X$ be distinct points. 
  Urysohn's lemma gives $f\in C(X)$ with $f(x)=1\neq 0=f(y)$. Then there exist a finite-dimensional C*-subalgebra $B$ and $g\in B$ with $\|f-g\|<\frac{1}{2}$. This implies that $g(x)\neq g(y)$.
  Since $B$ is finite-dimensional, Lemma \ref{lem:topologicalconditionoffinitedimensionalsubalgebra} below makes the set $\{ z \in X \mid \forall h \in B \colon h(x)=h(z) \}$ clopen. Hence $x$ and $y$ cannot share a connected component, and $X$ is totally disconnected. 

  Next we show (3)$\implies$(4), so let $X$ be totally disconnected. Distinct points $x,y\in X$ induce a clopen subset $C\subseteq X$ containing $x$ but not $y$. Now, the characteristic function of $C$ is continuous, and is therefore a projection, which clearly attains different values on $x$ and $y$. Thus the projections of $C(X)$ separate $X$. 
  It follows from the Stone--Weierstrass theorem~\cite[Theorem 3.4.14]{kadisonringrose:oa1} that the projections span an algebra $B$ that is dense in $C(X)$, hence $C(X)$ is generated by its projections.
  
  Finally, we prove (4)$\implies$(1), so assume that $C(X)$ is generated by its projections. Let $\D$ be the collection of all C*-subalgebras of $C(X)$ generated by a finite number of projections. By Lemma~\ref{lem:finiteprojectionsgeneratefinitedimsubalgebra}, $\D$ consists solely of finite-dimensional C*-subalgebras of $C(X)$. It is also directed: if $D_1,D_2\in\D$ where $D_1$ is generated by projections $p_1,\ldots,p_n$ and $D_2$ is generated by projections $q_1,\ldots,q_m$, then the C*-subalgebra $D$ generated by $p_1,\ldots,p_n,q_1,\ldots,q_m$ is clearly a member of $\D$ containing both $D_1$ and $D_2$. Since $\bigcup\D$ contains all projections of $C(X)$, and the latter C*-algebra is generated by its projections, it follows that $C(X)$ is the least upper bound of $\D$, hence $C(X)$ is approximately finite-dimensional.
\end{proof}

\begin{example}\label{ex:cantor}
  Let $X$ be the Cantor set. Then $C(X)$ is a separable commutative AF-algebra. Since there exists a continuous surjection $X\to[0,1]$, there is a C*-subalgebra of $C(X)$ that is $*$-isomorphic to $C([0,1])$ by Proposition~\ref{prop:dualitysubalgebrasandquotients}. This C*-subalgebra is not approximately finite-dimensional because $[0,1]$ is not totally disconnected. This pathological behaviour demonstrates why the notion of C*-subalgebras is sometimes replaced by that of hereditary C*-subalgebras, i.e., C*-subalgebras $B$ such that for each self-adjoint $b\in B$ and each self-adjoint $a$ in the ambient algebra the inequality $0\leq a\leq b$ implies $a\in B$, where $a\leq b$ if and only if $b-a=c^*c$ for some $c$ in the ambient algebra.
\end{example}

The following useful lemma explains the terminology `approximately' by linking the topology of a C*-algebra to approximating C*-subalgebras.

\begin{lemma}\label{lem:projectionsinAF}\cite[Proposition L.2.2]{weggeolsen:ktheory}
  Let $A$ be a C*-algebra and $\D$ a directed family of C*-subalgebras with $A=\overline{\bigcup\D}$. 
  For each $a\in A$ and $\varepsilon>0$, there exist $D\in\D$ and $x\in D$ satisfying $\|a-x\|<\varepsilon$. 
  If $a$ is a projection, then $x$ can be chosen to be a projection as well.
\end{lemma}

As far as we know, approximate finite-dimensionality of $A$ does not correspond to any nice order-theoretic properties of $\CC(A)$. We will need the following more subtle notion. In general, we will rely on the point-set topology of totally disconnected spaces, as covered \textit{e.g.}~in~\cite{gierzetal:domains}.

\begin{definition}\label{def:topologically scattered}
  A topological space is called \emph{scattered} if every nonempty closed subset has an isolated point.   
\end{definition}

Equivalently, a topological space $X$ is scattered if there is no continuous surjection $X \to [0,1]$~\cite[Theorem~8.5.4]{semadeni:banachspaces}.
Scattered topological spaces are always totally disconnected, so commutative C*-algebras with scattered spectrum are always approximately finite-dimensional. 

\begin{example}\label{ex:scattered}
  Any discrete topological space is scattered, and any finite discrete space is additionally compact Hausdorff, but there are more interesting examples.

  The one-point compactification of the natural numbers is scattered, as well as compact Hausdorff. This is homeomorphic to the subspace $\{ \frac{1}{n} \mid n \in \mathbb{N} \} \cup \{0\}$ of $\mathbb{R}$ under the usual Euclidean topology.

  More generally, any ordinal number $\alpha$ is scattered under the order topology. A basis for this topology is given by the intervals $\{ \delta \mid \beta < \delta < \gamma \}$ for ordinals $\beta,\gamma \leq \alpha$.
  If $\alpha$ is a limit ordinal, then $\alpha+1$ is furthermore compact Hausdorff \cite[Corollary 8.6.7]{semadeni:banachspaces}.
\end{example}

There is also an established notion of scatteredness in general C*-algebras $A$, which can be defined as follows.
\begin{definition}\cite{jensen:scattered}. 
A \emph{positive functional} on a C*-algebra $A$ is a continuous linear map $f \colon A\to\C$ satisfying $f(a^*a)\geq 0$. Positive functionals of unit norm are called \emph{states}, and form a convex set, whose extremal points are called \emph{pure}. A positive functional is called \emph{pure} if it is a positive multiple of a pure state.
$A$ is called a \emph{scattered} C*-algebra when each positive functional can be written as the countable sum of pure positive functionals, where the sum converges pointwise. 
\end{definition}
Theorem~\ref{thm:Kusuda} characterizes scattered C*-algebras completely. For now, let us mention that scattered C*-algebras are locally AF-algebras (cf.~\cite[Lemma 5.1]{lin:quasimultipliers}, where the author uses `AF-algebra' to mean 'locally AF-algebra'); it was only recently that a scattered C*-algebra that is not approximately finite-dimensional was found~\cite[Theorem 1.10]{bicekoszmider:approximateunits}.

\begin{example}\label{ex:compact operators}
  An operator $f \in B(H)$ on a Hilbert space $H$ is \emph{compact} when it is a limit of operators of finite rank.
  If $H$ is infinite-dimensional, the compact operators form a proper ideal $K(H)\subseteq B(H)$,
  and all self-adjoint elements of $K(H)$ have countable spectrum~\cite[Theorem VII.7.1]{conway:functionalanalysis}.
  It follows that the C*-algebra $K(H) + \C 1_H$ is scattered~\cite{huruya:spectral}.
\end{example}

The following theorem connects the notions of AF-algebras, scattered topological spaces, and scattered C*-algebras.

\begin{theorem}\label{thm:Kusuda}
  The following are equivalent for C*-algebras:
  \begin{enumerate}
    \item[(1)] $A$ is scattered;
    \item[(2)] each $C \in \CC(A)$ is approximately finite-dimensional;
    \item[(3)] each $C \in \CC(A)$ has totally disconnected spectrum;
    \item[(4)] each maximal $C \in \CC(A)$ has scattered spectrum;
    \item[(5)] no $C \in \CC(A)$ has spectrum $[0,1]$.
  \end{enumerate}
\end{theorem}
\begin{proof}  
  It follows from~\cite[Theorem 2.2]{kusuda:scattered} that (1) implies (2). The converse follows from \cite[Theorem 2.3]{kusuda:AF}. The equivalence between (2) and (3) is proven in Proposition~\ref{prop:commutativeAFistotallydisconnected}. By the same proposition, (2) implies (5). We show that (5) implies (4) by contraposition. Assume that some maximal commutative C*-subalgebra $M$ has non-scattered spectrum $X$. Then there is a continuous surjection $X \to [0,1]$, and it follows that $M$, and hence $A$, has a commutative C*-subalgebra whose spectrum is (homeomorphic to) $[0,1]$. Finally, we show that (4) implies (3).  Assume that all maximal commutative C*-subalgebras have scattered spectrum. Since by Zorn's lemma every commutative C*-subalgebra is contained in a maximal one, it follows from~\cite[Lemma 12.24]{fabianetal:functionalanalysis} that all commutative C*-subalgebras have scattered and hence totally disconnected spectrum.
 \end{proof}

\section{Algebraicity}\label{sec:algebraicity}

In this section we characterize C*-algebras $A$ for which $\CC(A)$ is algebraic. First recall what the latter notion means.
Consider elements $B,C$ of a dcpo $\CC$. The element $C$ could contain so much information that it is practically unobtainable. What does it mean for $B$ to approximate $C$ empirically? One answer is: whenever $C$ is the eventual observation of increasingly fine-grained experiments $D_i$, then all information in $B$ is already contained in a single one of the approximants $D_i$. 
More precisely: we say that $B$ is \emph{way below} $C$ and write $B\ll C$ if for each directed subset $\{D_i\}$ of $\CC$ the inequality $C\leq\bigvee D_i$ implies that $B \leq D_i$ for some $i$. Define:
\[
  \ddown C = \{B\in\CC \mid B\ll C\},
  \qquad
  \uup C=\{B\in\CC \mid C\ll B\}.
\]
With this interpretation, $C$ is empirically accessible precisely when $C \ll C$. Such elements are called \emph{compact}, and the subset they form is denoted by $\K(\CC)$.

\begin{definition}
  A dcpo is \emph{algebraic} when each element $C$ satisfies $C = \bigvee^{\up}\big( \K(\CC) \cap \down C \big)$, i.e., if every element is the supremum of the compact elements below it.
\end{definition}

We start by identifying the compact elements of $\CC(A)$. 
If $K$ is a closed subspace of a compact Hausdorff space $X$, define 
\[
  C_K = \{f\in C(X) \mid f~\text{is~constant~on}~K\},
\]
which is clearly a C*-subalgebra of $C(X)$. The following lemma gives a convenient way to construct directed subsets of $\CC(A)$.

\begin{lemma}\label{lem:directedfamilygeneratedbylocalneighbourhoods}
  Let $A$ be a C*-algebra, and $C\subseteq A$ a commutative C*-subalgebra with spectrum $X$. If $P\subseteq X$ is finite, then
  \[
   \D=\left\{\bigcap\nolimits_{p\in P}C_{\overline{U_p}} \;\;\big|\; U_p \text{ open neighbourhood of } p \right\}
  \]
  is a directed family in $\CC(A)$ with supremum $C$.
\end{lemma}
\begin{proof}
  If $U_p$ and $V_p$ are open neighbourhoods of $p$, then so is $U_p\cap V_p$.
  Moreover
  \[
    \overline{U_p\cap V_p} \subseteq \overline{U_p} \cap \overline{V_p} \subseteq \overline{U_p},
  \]
  and so $C_{\overline {U_p}} \subseteq C_{\overline{U_p\cap V_p}}$. 
  Similarly $C_{\overline {V_p}}\subseteq C_{\overline{U_p\cap V_p}}$. 
  Hence $\bigcap C_{\overline{U_p}}$ and $\bigcap C_{\overline {V_p}}$ are both contained in $\bigcap C_{\overline{U_p\cap V_p}}$, making $\D$ directed.
   
  Since Proposition \ref{prop:CAisdcpo} assures that $\bigvee\D$ is a C*-subalgebra of $C=C(X)$, it is closed under the operation $f\mapsto f^*$, and it contains the identity of $C(X)$, hence also all constant functions on $X$. Hence to show that $\bigvee\D=C(X)$, it suffices to show that $\bigvee\D$ separates all points of $X$ by the Stone--Weierstrass Theorem~\cite[Theorem 3.4.14]{kadisonringrose:oa1}.
  Thus let $x$ and $y$ be distinct points in $X$; we will show that $f(x) \neq f(y)$ for some $f\in\bigvee\D$ by distinguishing two cases.
  For the first case, suppose $x,y\in P$. 
  Since $P$ is finite, it is closed, as is $P\setminus\{x\}$.
  Hence $\{x\}$ and $P\setminus\{x\}$ are disjoint closed subsets in $X$, and since $X$ is compact Hausdorff, there are open subsets $U$ and $V$ containing $x$ and $P\setminus\{x\}$, respectively, with disjoint closures.
  Because $U$ is an open neighbourhood of $x$ and $V$ is an open neighbourhood of $p$ for each $p\in P\setminus\{x\}$, it follows that $C_{\overline U}\cap C_{\overline V}$ is in $\D$. 
  But Urysohn's lemma provides a function $f\in C(X)$ satisfying $f(\overline U)=\{0\}$ and $f(\overline V)=\{1\}$. 
  Hence $f$ is constant on $\overline U$ and on $\overline V$, so $f\in C_{\overline U}\cap C_{\overline V}$. 
  Since $y\in P\setminus\{x\}\subseteq\overline V$, we find $f(x)=0\neq 1=f(y)$.
  
  For the second case, suppose $x\not\in P$, and proceed similarly.
  Regardless of whether $y\in P$ or not, $\{x\}$ and $P\cup\{y\}$ are disjoint closed subsets, hence there are open sets $U$ and $V$ containing $\{x\}$ and $P\cup\{y\}$, respectively, with disjoint closures.
  Since $V$ is an open neighbourhood of $p$ for each $p\in P$, we find that $C_{\overline V}$ is in the family. 
  Again Urysohn's lemma provides a function $f\in C(X)$ satisfying $f(\overline U)=\{0\}$ and $f(\overline V)=\{1\}$, and since $f$ is constant on $\overline V$, we find $f\in C_{\overline V}$, and again $f(x)\neq f(y)$.
\end{proof}

We can now identify the compact elements of $\CC(A)$ as the finite-dimensional ones.
 
\begin{proposition}\label{prop:finitedimiscompact}
  Let $A$ be a C*-algebra. Then $C\in\CC(A)$ is compact if and only if it is finite-dimensional.
\end{proposition}
\begin{proof}
  Suppose $C$ is compact, and write $X$ for its spectrum.
  Let $x\in X$ and consider 
  \[
    \D = \{C_{\overline U} \mid U \text{ is an open neighbourhood of } x \}. 
  \]
  It follows from Lemma \ref{lem:directedfamilygeneratedbylocalneighbourhoods} that $\D$ is directed and $C(X)=\bigvee\D$. 
  Because $C$ is compact, it must equal some element $C_{\overline{U}}$ of $\D$.
  Since $C(X)$ separates all points of $X$, so must $C_{\overline{U}}$. 
  But as each $f\in C_{\overline U}$ is constant on $\overline U$, this can only happen when $\overline{U}$ is a singleton $\{ x \}$.
  This implies $\{x\}=U$, so $\{x\}$ is open. 
  Since $x\in X$ was arbitrary, $X$ must be discrete. 
  Being compact, it must therefore be finite.
  Hence $C$ is finite-dimensional.
  
  Conversely, assume that $C$ has a finite dimension $n$. By Theorem \ref{thm:gelfandduality}, $C$ is isomorphic to $C(X)$, where $X$ is a discrete $n$-point space, which is clearly spanned by the characteristic functions on the singleton sets. Hence $C$ is generated by a finite set $\{p_1,\ldots,p_n\}$ of \emph{mutually orthogonal} projections in the sense that $p_i p_j=0$ for $i \neq j$.
  Let $\D\subseteq\CC(A)$ be a directed family satisfying $C\subseteq\bigvee\D$, and let $X$ be the Gelfand spectrum of $\bigvee\D$. Fix $i\in\{1,\ldots,n\}$. 
  Since the projection $p_i$ is contained in $\bigvee\D$, using Lemma~\ref{lem:projectionsinAF} there is some $D\in\D$ and some projection $p\in D$ such that $\|p_i-p\|<\frac{1}{2}$.
  Since projections $p \colon X \to \C$ can only take the value 0 or 1, $p \neq p_i$ implies $\|p_i-p\|=1$, so we must have $p=p_i$. 
  Hence there are $D_{1},\ldots,D_{n}\in\D$ such that $p_i\in D_i$ for each $i\in\{1,\ldots,n\}$. 
  Since $\D$ is directed, there must be some $D\in\D$ with $D_1,\ldots,D_n\subseteq D$. So $p_1,\ldots,p_n \in D$, which implies that $C \subseteq D$. 
  We conclude that $C$ is compact.
\end{proof}

This leads to the following characterization of algebraicity of $\CC(A)$.

\begin{theorem}\label{thm:algebraic}
  A C*-algebra $A$ is scattered if and only if $\CC(A)$ is algebraic.
\end{theorem}
\begin{proof}
  By Proposition~\ref{prop:finitedimiscompact} and Proposition~\ref{prop:commutativeAFistotallydisconnected}, the dcpo $\CC(A)$ is algebraic if and only if each $C\in\CC(A)$ is approximately finite-dimensional. By Theorem~\ref{thm:Kusuda}, this is equivalent with scatteredness of $A$.
\end{proof}

\section{Continuity}\label{sec:continuity}

In this section we characterize C*-algebras $A$ for which $\CC(A)$ is continuous. 

\begin{definition}
  A dcpo is \emph{continuous} when each element satisfies $C= \bigvee^{\up} \ddown C$.
\end{definition}

We start with two lemmas that govern the equivalence relation $\sim_B$ on a compact Hausdorff space $X$ defined by $x \sim_B y$ if and only if $b(x)=b(y)$ for each element $b$ of a C*-subalgebra $B \subseteq C(X)$.

\begin{lemma}\label{lem:topologicalconditionoffinitedimensionalsubalgebra}
  For a compact Hausdorff space $X$ and C*-subalgebra $B\subseteq C(X)$:
  \begin{enumerate}
  	\item[(1)] each equivalence class $[x]_B$ is a closed subset of $X$.
    \item[(2)] $B$ is finite-dimensional if and only if $[x]_B\subseteq X$ is open for each $x\in X$;
    \item[(3)] if $X$ is connected, $B$ is the (one-dimensional) subalgebra of all constant functions on $X$ if and only if $[x]_B$ is open for some $x\in X$;
    \item[(4)] if $B$ is infinite-dimensional, there are $x\in X$ and $p\in [x]_B$ such that $B\nsubseteq C_{\overline{U}}$ on each open neighbourhood $U\subseteq X$ of $p$. If $X$ is connected, this holds for all $x\in X$.
  \end{enumerate}
\end{lemma}
\begin{proof*}
  Fix $X$ and $B$.
  \begin{enumerate}
  	\item[(1)] The proof of Proposition~\ref{prop:dualitysubalgebrasandquotients} shows that the quotient $X/{\sim_B}$ is compact Hausdorff, hence its points are closed. 
  	If $q$ denotes the quotient map, then $[x]_B$ as subset of $X$ is equal to the preimage under the continuous map $q$ of $[x]_B$ as point of $X/{\sim_B}$. Thus $[x]_B$ is closed.
    \item[(2)]  Let $q \colon X\to X/{\sim_B}$ be the quotient map. 
    By definition of the quotient topology, $V\subseteq X/{\sim_B}$ is open if and only if its preimage $q^{-1}[V]$ is open in $X$. 
    We can regard $[x]_B$ both as a subset of $X$ and as a point in $X/{\sim_B}$. 
    Since $[x]_B=q^{-1}[\{[x]_B\}]$, we find that $\{[x]_B\}$ is open in $X/{\sim_B}$ if and only if $[x]_B$ is open in $X$. 
    Hence $X/{\sim_B}$ is discrete if and only if $[x]_B$ is open in $X$ for each $x\in X$. 
    Now $X/{\sim_B}$ is compact, being a continuous image of a compact space.
    It is also Hausdorff by Proposition~\ref{prop:dualitysubalgebrasandquotients}.
    Hence $X/{\sim_B}$ is discrete if and only if it is finite. 
    Thus each $[x]_B$ is open in $X$ if and only if $B$ is finite-dimensional.
    
    \item[(3)] An equivalence class $[x]_B$ is always closed in $X$ by (1). Assume that it is also open.
    By connectedness $X=[x]_B$, so $f(y)=f(x)$ for each $f\in B$ and each $y\in X$.
    Hence $B$ is the algebra of all constant functions on $X$, and since this algebra is spanned by the function $x\mapsto 1$, it follows that $B$ is one dimensional. 
    
    Conversely, if $B$ is the one-dimensional subalgebra of all constant function on $X$, then for each $f\in B$ there is some $\lambda\in\C$ such that $f(x)=\lambda$ for each $x\in X$. Hence $f(x)=f(y)$ for each $x,y\in X$, whence for each $x\in X$ we have $[x]_B=X$, which is clearly open.
    
    \item[(4)] Assume that $B$ is infinite-dimensional. 
    By (2) there must be some $x\in X$ such that $[x]_B$ is not open.
    Hence there must be a point $p\in[x]_B$ such that $U\nsubseteq [x]_B$ for each open neighbourhood $U$ of $p$. If $X$ is connected, (3) implies that $[x]_B$ is not open for any $x\in X$, so $p$ can be chosen as an element of $[x]_B$ for each $x\in X$. 
    In both cases, we have $U\nsubseteq[x]_B$ for any open neighbourhood $U$ of $p$, hence there is $q\in U$ such that $q\notin[x]_B$. 
    We have $y\in[x]_B$ if and only if $f(x)=f(y)$ for each $f\in B$. So $p\in[x]_B$, and $q\notin[x]_B$ implies the existence of some $f\in B$ such that $f(p)\neq f(q)$. That is, there is some $f\in B$ such that $f$ is not constant on $U$, so $f$ is certainly not constant on $\overline{U}$. We conclude that for each open neighbourhood $U$ of $p$ there is an $f\in B$ such that $f\notin C_{\overline U}$, so $B\nsubseteq C_{\overline U}$ for each open neighbourhood $U$ of $p$.
    {\hspace*{\fill}{\usebox{\proofbox}}\endlist}
  \end{enumerate}
\end{proof*}

We can now characterize the way-below relation on $\CC(A)$ in operator-algebraic terms.

\begin{proposition}\label{prop:waybelowinCA}
  The following are equivalent for a C*-algebra $A$ and $B,C\in\CC(A)$:
  \begin{enumerate}
    \item[(1)] $B\ll C$;
    \item[(2)] $B\in\K(\CC)$ and $B\subseteq C$;
    \item[(3)] $B$ is finite-dimensional and $B\subseteq C$.
  \end{enumerate}
\end{proposition}
\begin{proof}
  By Proposition \ref{prop:finitedimiscompact}, $B$ is finite-dimensional if and only if $B$ is compact, which proves the equivalence between (2) and (3). 
  It is almost trivial that (2) implies (1) by unfolding definitions.
  For (1) $\Rightarrow$ (3), assume $B\ll C$, which implies $B\subseteq C$. For a contradiction, assume that $B$ is infinite-dimensional. 
  Without loss of generality we may assume that $C=C(X)$ for the spectrum $X$ of $C$.
  Lemma~\ref{lem:topologicalconditionoffinitedimensionalsubalgebra} gives $p\in X$ with $B\nsubseteq C_{\overline U}$ for each open neighbourhood $U \subseteq X$ of $p$. Consider the family
  \[
    \{C_{\overline U} \mid U\text{ open neighbourhood of }p\}\text{.}
  \]
  By Lemma~\ref{lem:directedfamilygeneratedbylocalneighbourhoods}, this is a directed family in $\CC(A)$ with supremum $C(X)$. 
  However, $B$ is not contained in any member of the family, and so cannot be way below $C=C(X)$.
\end{proof}

This leads to the following characterization of continuity of $\CC(A)$.

\begin{theorem}\label{thm:continuous}
  A C*-algebra $A$ is scattered if and only if $\CC(A)$ is continuous.
\end{theorem}
\begin{proof}
  Let $C\in\CC(A)$. It follows from Proposition~\ref{prop:waybelowinCA} that $\ddown C=\K(\CC)\cap\down C$, whence $C=\bigvee\K(\CC)\cap\down C$ if and only if $C=\bigvee\ddown C$. Thus continuity and algebraicity coincide. The statement now follows from Theorem~\ref{thm:algebraic}.
\end{proof}

\section{Meet-continuity}\label{sec:meetcontinuity}

In this section we characterize C*-algebras $A$ for which $\CC(A)$ is meet-continuous.
  
\begin{definition}
  A dcpo $\CC$ is \emph{meet-continuous} when it is a meet-semilattice, and
  \begin{equation}\label{eq:distributivityifAfinitedimensional}
    C \wedge \bigvee\D = \bigvee_{D\in\D} C \wedge D
  \end{equation}
  for each element $C$ and directed subset $\D$ of $\CC$.
\end{definition}

A \emph{closed equivalence relation} on a topological space $X$ is a reflexive, symmetric, and transitive relation $R$ on $X$ that is closed as a subset of $X \times X$ in the product topology.

\begin{lemma}\label{lem:equivalencerelations}
  There is a dual equivalence between the poset of closed equivalence relations on a compact Hausdorff space $X$ and $\mathcal{C}(C(X))$, that sends $R$ to 
  \[
    C_R = \{ f \in C(X) \mid \forall x,y \in X \colon (x,y) \in R \Rightarrow f(x)=f(y) \} \text{.}
  \]
\end{lemma}
\begin{proof}
  The map $R \mapsto C_R$ is essentially bijective because any $C \in \mathcal{C}(C(X))$ corresponds to a quotient $X \twoheadrightarrow \Spec(C)$ of compact Hausdorff spaces, which in turn corresponds to a closed equivalence relation $\sim$ on $X$ by $\Spec(C) = X \mathop{/} \mathrm{\sim}$. Clearly $R \subseteq S$ if and only if $C_R \supseteq C_S$.
\end{proof}

Closed equivalence relations on $X$ form a complete lattice under reverse inclusion, which follows from the observation that the intersection of a family of equivalence relations is an equivalence relation, and the intersection of a family of closed subsets is closed, hence the infimum of a family of closed equivalence relations is simply given by intersection. The supremum is harder to describe, and can in general only be given as
\[
  \bigvee R_n = \bigcap \left\{ S \subseteq X^2 \text{ closed equivalence relation} \;\big|\; \bigcup R_n \subseteq S\right\}\text{.}
\]
Recall that composition of relations $R$ and $S$ on $X$ is defined by \[R\circ S=\{(x,z)\in X^2 \mid \exists y\in X\colon (x,y)\in S, (y,z)\in R\}\text.\] 
In the special case where $R \circ S = S \circ R$ we have $R \vee S = R \circ S$~\cite[Proposition~6.9]{ellisellis:equivalencerelations}. 
By Lemma~\ref{lem:equivalencerelations}, meet-continuity of $\CC(C(X))$ comes down to the question whether $R \vee \bigcap S_n \supseteq \bigcap R \vee S_n$ for closed equivalence relations $R$ and $S_1 \supseteq S_2 \supseteq S_3 \supseteq \cdots$ on $X$.
Notice that $C_K$, for a C*-algebra $C=C(X)$ and closed subset $K \subseteq X$, is a special case of $C_R$ of Lemma~\ref{lem:equivalencerelations} for the closed equivalence relation $R=\{(x,x) \mid x \in X\} \cup K^2 \subseteq X^2$. 
In general $C_R = \bigcap_{x \in X} C_{[x]_R}$.
It will always be clear from the context which of the two is meant.

\begin{proposition}\label{prop:equivalencerelations}
  There are closed equivalence relations $R$ and $S_1 \supseteq S_2 \supseteq S_3 \supseteq \cdots$ on $[0,1]$ such that
  $R \vee \bigcap S_n \neq \bigcap R \vee S_n$.
\end{proposition}
\begin{proof}
  We will construct $R$ and $S_n$ similar to the Cantor set; $R$ by keeping the closed middle thirds, and $S_n$ by keeping the closed first and last thirds. The definitions of $S_n$ below can be illustrated as follows:
  \[
    S_1 = 
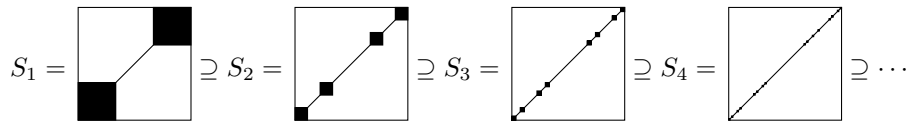
\begin{aligned}\begin{tikzpicture}[scale=1.5]
      \draw (0,0) rectangle (1,1); \draw (0,0) to (1,1);
    \draw[fill] (0/3,0/3) rectangle (1/3,1/3);
    \draw[fill] (2/3,2/3) rectangle (3/3,3/3);
    \end{tikzpicture}\end{aligned}
    \supseteq 
    S_2 = \begin{aligned}\begin{tikzpicture}[scale=1.5]
      \draw (0,0) rectangle (1,1); \draw (0,0) to (1,1);
    \draw[fill] (0/9,0/9) rectangle (1/9,1/9);
    \draw[fill] (2/9,2/9) rectangle (3/9,3/9);
    \draw[fill] (6/9,6/9) rectangle (7/9,7/9);
    \draw[fill] (8/9,8/9) rectangle (9/9,9/9);
    \end{tikzpicture}\end{aligned}
    \supseteq
    S_3 = \begin{aligned}\begin{tikzpicture}[scale=1.5]
      \draw (0,0) rectangle (1,1); \draw (0,0) to (1,1);
    \draw[fill] (0/27,0/27) rectangle (1/27,1/27);
    \draw[fill] (2/27,2/27) rectangle (3/27,3/27);
    \draw[fill] (6/27,6/27) rectangle (7/27,7/27);
    \draw[fill] (8/27,8/27) rectangle (9/27,9/27);
    \draw[fill] (18/27,18/27) rectangle (19/27,19/27);
    \draw[fill] (20/27,20/27) rectangle (21/27,21/27);
    \draw[fill] (24/27,24/27) rectangle (25/27,25/27);
    \draw[fill] (26/27,26/27) rectangle (27/27,27/27);
    \end{tikzpicture}\end{aligned}
    \supseteq
    S_4 = \begin{aligned}\begin{tikzpicture}[scale=1.5]
      \draw (0,0) rectangle (1,1); \draw (0,0) to (1,1);
    \draw[fill] (0/81,0/81) rectangle (1/81,1/81);
    \draw[fill] (2/81,2/81) rectangle (3/81,3/81);
    \draw[fill] (6/81,6/81) rectangle (7/81,7/81);
    \draw[fill] (8/81,8/81) rectangle (9/81,9/81);
    \draw[fill] (18/81,18/81) rectangle (19/81,19/81);
    \draw[fill] (20/81,20/81) rectangle (21/81,21/81);
    \draw[fill] (24/81,24/81) rectangle (25/81,25/81);
    \draw[fill] (26/81,26/81) rectangle (27/81,27/81);
    \draw[fill] (54/81,54/81) rectangle (55/81,55/81);
    \draw[fill] (56/81,56/81) rectangle (57/81,57/81);
    \draw[fill] (60/81,60/81) rectangle (61/81,61/81);
    \draw[fill] (64/81,64/81) rectangle (65/81,65/81);
    \draw[fill] (72/81,72/81) rectangle (73/81,73/81);
    \draw[fill] (74/81,74/81) rectangle (75/81,75/81);
    \draw[fill] (78/81,78/81) rectangle (79/81,79/81);
    \draw[fill] (80/81,80/81) rectangle (81/81,81/81);
    \end{tikzpicture}\end{aligned}
    \supseteq \cdots
  \]
  Similarly, $R$ can be drawn as follows:
  \[
  R = 
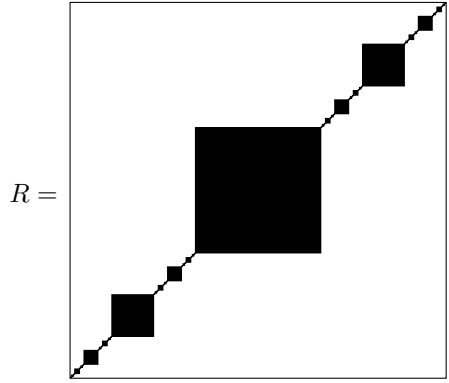
\begin{aligned}\begin{tikzpicture}[scale=5]
    \draw (0,0) rectangle (1,1); \draw (0,0) to (1,1);
    \foreach \n in {1,...,6}
    {
      \pgfmathsetmacro{\N}{3^(\n-1)-1};
      \pgfmathsetmacro{\d}{3^\n};
        \foreach \k in {0,...,\N}
        {
          \pgfmathsetmacro{\l}{(3*\k+1) / \d};
          \pgfmathsetmacro{\r}{(3*\k+2) / \d};
          \draw[fill] (\l,\l) rectangle (\r,\r);
        }   
      }
  \end{tikzpicture}\end{aligned}
  \]
  Formally, let $\epsilon$ denote the empty string. Inductively define numbers $a_\sigma,b_\sigma,c_\sigma,d_\sigma \in [0,1]$ indexed by finite strings $\sigma$ of zeroes and ones:
  \begin{align*}
  a_\epsilon & = 0, &
  b_\epsilon & = 1/3,&
  c_\epsilon & = 2/3,&
  d_\epsilon & = 1, \\
  a_{\sigma 0} & = a_\sigma,&
  b_{\sigma 0} & = a_\sigma + \tfrac{1}{3}(b_\sigma-a_\sigma),&
  c_{\sigma 0} & = b_\sigma - \tfrac{1}{3}(b_\sigma-a_\sigma),&
  d_{\sigma 0} & = b_\sigma, \\
  a_{\sigma 1} & = c_\sigma,&
  b_{\sigma 1} & = c_\sigma + \tfrac{1}{3}(d_\sigma-c_\sigma),&
  c_{\sigma 1} & = d_\sigma - \tfrac{1}{3}(d_\sigma-c_\sigma),&
  d_{\sigma 1} & = d_\sigma,
  \end{align*} 
  where $\sigma\in\{0,1\}^*$, where the star is the Kleene star.
  Write $\Delta = \{ (x,x) \mid x \in [0,1] \}$ for the diagonal, and define
  \begin{align*}
    R & = \Delta \cup \bigcup_{\sigma \in \{0,1\}^*} [b_\sigma, c_\sigma]^2\text{,} \\
    S_n & = \Delta \cup \bigcup_{\sigma \in \{0,1\}^n} [a_\sigma,d_\sigma]^2\text{.}
  \end{align*}
  The $S_n$ are certainly closed equivalence relations, and $\bigcap S_n = \Delta$.

  Clearly $R$ is reflexive and symmetric. It is also transitive: if $(x,y)$ and $(y,z)$ in $R$ have $x \neq y \neq z$, then $(x,y) \in [b_\sigma,c_\sigma]^2$ and $(y,z) \in [b_\tau,c_\tau]^2$ for some $\sigma,\tau \in \{0,1\}^*$; but if $y \in [b_\sigma,c_\sigma] \cap [b_\tau,c_\tau]$ then $\sigma=\tau$, so $(x,z) \in R$.
  The set $R \subseteq [0,1]^2$ is also closed: if $(x_n,y_n) \in R$ is a sequence that converges in $[0,1]^2$, then either it eventually stays in one block $[b_\sigma,c_\sigma]^2$, or it converges to a point on the diagonal.

  In total we see that $R \vee \bigcap S_n = R \vee \Delta = R$.
  But we now prove that $R \vee S_n = [0,1]^2$ for any $n$, so $\bigcap R \vee S_n = [0,1]^2$, and hence $R \vee \bigcap S_n \neq \bigcap R \vee S_n$. 
  By induction it suffices to show $R \vee S_1 = [0,1]^2$ and $S_n \subseteq R \vee S_{n+1}$. For the latter it suffices to show $(a_\sigma,d_\sigma) \in R \vee S_{n+1}$ for $\sigma \in \{0,1\}^n$, which follows from transitivity:
  \[
    a_\sigma 
    = a_{\sigma 0}
    \,S_{n+1}\,
    b_{\sigma 0}
    \,R\,
    c_{\sigma 0}
    \,S_{n+1}\,
    d_{\sigma 0}
    = b_\sigma
    \,R\,
    c_\sigma
    = a_{\sigma 1}
    \,S_{n+1}\,
    b_{\sigma 1}
    \,R\,
    c_{\sigma 1}
    \,S_{n+1}\,
    d_{\sigma 1}
    = d_\sigma.
  \]
  Similarly $R \vee S_1 = [0,1]^2$.
\end{proof}

\begin{theorem}
  A C*-algebra $A$ is scattered if and only if $\CC(A)$ is meet-continuous.
\end{theorem}
\begin{proof}
  If $A$ is not scattered, there is an element of $\CC(A)$ that is $*$-isomorphic to $C([0,1])$ by Theorem~\ref{thm:Kusuda}. Therefore we may assume without loss of generality that $A=C([0,1])$. But it now follows from Lemma~\ref{lem:equivalencerelations} and Proposition~\ref{prop:equivalencerelations} that $\CC(A)$ is not meet-continuous.

  If $A$ is scattered, then $\CC(A)$ is continuous by Theorem~\ref{thm:continuous}. But $\CC(A)$ is also a complete semilattice, because the intersection $\bigcap C_i$ of a family of commutative C*-subalgebras $C_i$ of $A$ is again a commutative C*-subalgebra. And continuous dcpos that are also semilattices are meet-continuous~\cite[Proposition I-1.8]{gierzetal:domains}.
\end{proof}

We can give another characterization of meet-continuity of $\CC(A)$, namely in order-theoretic terms of injective $*$-homomorphisms $f:B \to A$. In this case the upper adjoint of $\CC(f) \colon \CC(A)\to\CC(B)$, \ie\ the monotone map $\CC(f)_*:\CC(B)\to\CC(A)$ satisfying 
\[
	\CC(f)(C)\subseteq D\iff C\subseteq\CC(f)_*(D)
\]
exists, and is given by $D\mapsto f^{-1}[D]$, which clearly preserves inclusions. To see that this map is well defined, we have to show that $f^{-1}[D]$ is a commutative C*-subalgebra of $A$ if $D$ is a commutative C*-subalgebra of $B$. Clearly $f^{-1}[D]$ is a *-subalgebra that topologically closed (as $*$-homomorphisms are continuous). To show that it is commutative, let $x,y\in f^{-1}[D]$. Then \[f(xy-yx)=f(x)f(y)-f(y)f(x)=0,\] since $D$ is commutative. By injectivity of $f$ it follows that $xy-yx=0$, whence $f^{-1}[D]$ is commutative. To see that $\CC(f)_*$ is indeed the upper adjoint of $\CC(f)$, we recall that the latter is given by $C\mapsto f[C]$, and it is well known that \[f[C]\subseteq D\iff C\subseteq f^{-1}[D].\] 

\begin{theorem}\label{thm:meetcontinuity:maps}
  The dcpo $\CC(A)$ of a C*-algebra $A$ is meet-continuous if and only if $\CC(f)_*$ is Scott-continuous for any injective $*$-homomorphism $f \colon B \to A$.
\end{theorem}
\begin{proof}
  Suppose $\CC(f)_*$ is Scott-continuous for each injective $*$-homomorphisms $f \colon B \to A$.
  Let $\D$ be a directed family in $\CC(A)$, and let $C \in \CC(A)$. Write $i \colon C \to A$ for the inclusion.
  Then $\CC(i)$ has an upper adjoint $\CC(i)_* \colon \CC(B) \to \CC(A)$ given by $D\mapsto i^{-1}[D]$ so
  $\CC(i)_*(D) = i^{-1}[D] = C \cap D$.
  Because $\CC(i)_*$ is Scott-continuous, 
  \[
    C \cap \bigvee\D = \CC(i)_*\left(\bigvee\D \right) = \bigvee_{D\in\D} \CC(i)_*[D] = \bigvee_{D\in\D} C \cap D.
  \]
  Hence $\CC(A)$ is meet-continuous.

  Now assume $\CC(A)$ is meet-continuous and let $f \colon B\to A$ be an injective $*$-homomorphism.
  Write $\CC(f)_* \colon \CC(B) \to \CC(A)$ for the upper adjoint, and let $\D$ be a directed family in $\CC(B)$.
  Then $\CC(\varphi)_*(D) \subseteq \CC(f)_*\left(\bigvee\D \right)$ for each $D\in\D$, 
  and hence
  \[
    \bigvee_{D\in\D} f^{-1}[D]
    \;=\; \bigvee_{D\in\D} \CC(f)_*(D)
    \;\subseteq\; \CC(f)_* \left(\bigvee\D\right) 
    \;=\; f^{-1}\left[\bigvee\D \right].
  \]
  To show that this inclusion is in fact an equality, let  $x\in f^{-1}[\bigvee\D]$ be self-adjoint. Then the smallest C*-algebra $C=C^*(x)$ containing $x$ is commutative, hence an element of $\CC(B)$. 
 It follows from Lemma \ref{lem:imageofgeneratedsubalgebra} that $f[C]=C^*(f(x))$. 
   Since $f(x)\in\bigvee\D$, hence $C^*(f(x))\subseteq\bigvee\D$, it follows that $f[C]\subseteq\bigvee\D$. 
  Meet-continuity of $\CC(A)$ now shows $f[C] = \bigvee_{D\in\D} f[C] \cap D$.
  Being a C*-subalgebra, $f[B]$ is closed in $A$, so that the injection $f$ restricts to a $*$-isomorphism and hence a homeomorphism $B \to f[B]$.
  Observe that $f^{-1}[\overline{S}] = \overline{f^{-1}[S]}$ for $S \subseteq f[B]$.
  Hence
  \begin{align*}
    C 
    & = f^{-1}[f[C]]
      = f^{-1} \Big[ \bigvee\nolimits_{D\in\D} f[C]\cap D \Big]
      = f^{-1} \Big[ \overline{\bigcup\nolimits_{D\in\D} f[C]\cap D} \Big] \\
    & = \overline{f^{-1} \Big[ \bigcup\nolimits_{D\in\D} f[C]\cap D \Big]}
      = \overline{\bigcup\nolimits_{D\in\D} f^{-1}\big[f[C]\cap D\big]}  \\
    & \subseteq\; \overline{\bigcup\nolimits_{D\in\D} f^{-1}[D] }
      = \bigvee\nolimits_{D\in\D} f^{-1}[D] \text{.}
  \end{align*}
  As $x\in C$, it follows that $x\in\bigvee_{D\in\D} f^{-1}[D]$. 
  Finally, since $f^{-1}[\bigvee\D]$ is a C*-subalgebra of $A$, any of its elements is a linear combination of self-adjoint elements in $f^{-1}[\bigvee\D]$, hence $f^{-1}[\bigvee\D]\subseteq \bigvee_{D\in\D}f^{-1}[D]$.
  We conclude that
  $\CC(f)_* \left(\bigvee\D \right) = \bigvee_{D\in\D} \CC(f)_*(D)$. 
\end{proof}

\section{Atomicity}\label{sec:atomicity}

In this section we characterize the C*-algebras $A$ for which $\CC(A)$ is atomistic.

\begin{definition}
  Let $\CC$ be a partially ordered set with least element $0$. An \emph{atom} in $\CC$ is a minimal non-zero element. A partially ordered set is called \emph{atomistic} if each element is the least upper bound of some collection of atoms. 
\end{definition}

We begin by identifying the atoms in $\CC(A)$. Recall that $C^*(S)\subseteq A$ denotes the C*-subalgebra of $A$ generated by a subset $S$ of $A$, i.e., the smallest C*-subalgebra containing $S \subseteq A$. By Lemma \ref{lem:finiteprojectionsgeneratefinitedimsubalgebra}, $C^*(p)$ is just the linear span $\Span \{ p,1-p\}$ for projections $p^2=p^*=p \in A$; this is two-dimensional unless $p$ is \emph{trivial}, \ie 0 or 1, in which case it collapses to the least element $\C 1$ of $\CC(A)$.
The next lemma is simple and known~\cite{hamhalter:ordered}, but for completeness we include a proof.

\begin{lemma}\label{lem:atomsinCA}
  Let $A$ be a C*-algebra. Then $C$ is an atom in $\CC(A)$ if and only if it is generated by a nontrivial projection.
\end{lemma}
\begin{proof}
  Clearly two-dimensional $C$ are atoms in $\CC(A)$. 
  Conversely, assume that $C$ is an atom of $\CC(A)$. 
  By Theorem~\ref{thm:gelfandduality}, $C \simeq C(X)$ for a compact Hausdorff space $X$.
  If $X$ contains three distinct point $x,y,z$, then $C(X)$ contains a proper subalgebra $\{f\in C(X) \mid f(x)=f(y)\}$ with dimension at least two, which contradicts atomicity of $C$.
  Hence $X$ must contain exactly two points $x$ and $y$. Using the $*$-isomorphism between $C$ and $C(X)$, let $p\in C$ be the element corresponding to the element of $C(X)$ given by $x \mapsto 1$ and $y \mapsto 0$ for $y \neq x$.
  It follows that $C=\Span\{p,1-p\}$. 
\end{proof}

To characterize atomicity we will need two auxiliary results. The first deals with least upper bounds of subalgebras in terms of generators.

\begin{lemma}\label{lem:joinindownC}
  Let $A$ be a C*-algebra and $C\in\CC(A)$. If $\{S_i\}_{i\in I}$ is a family of subsets of $C$, then each $C^*(S_i)$ is in $\CC(A)$, and $C^*\left( \bigcup_{i \in I} S_i\right)=\bigvee_{i\in I}C^*(S_i)$.  
\end{lemma}
\begin{proof}
  For any $i \in I$, clearly $C^*(S_i)$ is a commutative C*-subalgebra of $A$, and hence an element of $\CC(A)$.

  Writing $S=\bigcup_{i \in I} S_i$, we have $S_j\subseteq C^*(S)$, and so $C^*(S_j)\subseteq C^*(C^*(S))=C^*(S)$. 
  Therefore, $\bigvee_{i\in I}C^*(S_i)$ is contained in $C^*(S)$. 
  For the inclusion in the other direction, notice that clearly $S\subseteq\bigvee_{i\in I}C^*(S_i)$, whence 
  \[
    C^*(S) \subseteq C^*\left(\bigvee_{i\in I}C^*(S_i)\right) = \bigvee_{i\in I}C^*(S_i).
  \]
  This finishes the proof.
\end{proof}

The second auxiliary result deals with subalgebras generated by projections. It shows that projections are the building blocks for C*-algebras $A$ whose dcpos $\CC(A)$ are atomistic. This explains why mere approximate finite-dimensionality is not good enough to characterize algebraicity and/or continuity. See also Section~\ref{sec:projections} below.

\begin{proposition}\label{prop:elementsofCgeneratedbyprojections}
  For a C*-algebra $A$, a C*-subalgebra $C$ of $A$ is the least upper bound of a collection of atoms of $\CC(A)$ if and only if it is generated by projections.
\end{proposition}
\begin{proof}
  Let $C\in\CC(A)$. For any collection $P$ of projections in $A$, of course $P=\bigcup_{p\in P}\{p\}$, so Lemma \ref{lem:joinindownC} guarantees $C^*(P)=\bigvee_{p\in P}C^*(p)$.  
  Hence $C$ is generated by projections if and only if $C=C^*(P)$ for some collection $P$ of projections in $C$ if and only if $C=\bigvee_{p\in P}C^*(p)$ for some collection of projections in  $C$. It now follows from Lemma~\ref{lem:atomsinCA} that $C$ is generated by projections if and only if it is the least upper bound of atoms in $\CC(A)$.
\end{proof}

This leads to the following characterization of atomicity of $\CC(A)$.

\begin{theorem}\label{thm:CAatomistic}
  A C*-algebra $A$ is scattered if and only if $\CC(A)$ is atomistic.  
\end{theorem}
\begin{proof}
  By Theorem~\ref{thm:algebraic} it suffices to prove that $\CC(A)$ is algebraic if and only if it is atomistic.
  Assume that $\CC(A)$ is algebraic and let $C\in\CC(A)$. 
  If $C=\C 1$, then $C$ is the least upper bound of the empty set, which is a subset of the set of atoms. 
  Otherwise, it follows from Proposition~\ref{prop:finitedimiscompact} that $C$ is the least upper bound of all its finite-dimensional C*-subalgebras. 
  Since every finite-dimensional C*-algebra is generated by a finite set of projections, it follows from Proposition~\ref{prop:elementsofCgeneratedbyprojections} that each element $D\in K(\CC(A))\cap\down C$ can be written as the least upper bound of atoms in $\CC(A)$. 
  Hence $C$ is a least upper bound of atoms, so $\CC(A)$ is atomistic.
  
  Conversely, assume $\CC(A)$ is atomistic and let $C\in\CC(A)$. 
  Because $C$ being finite-dimensional implies that it is a least upper bound of $K(\CC(A))\cap\down C$, we may assume that $C$ is infinite-dimensional. 
  By Lemma~\ref{lem:atomsinCA}, $C=\bigvee_{p\in P}C^*(p)$ for some collection $P$ of projections in $A$. As we must have $P\subseteq C$, all projections in $P$ commute. 
  We may replace $P$ by the set of all projections of $C$, which we will denote by $P$ as well; then we still have $C=\bigvee_{p\in P}C^*(p)$. 
  Write $\F$ for the collection of all finite subsets of $P$, and consider the family 
  $\{C^*(F) \in \CC(A) \mid F\in\F\}$.
  If $F\in\F$, then $C^*(F)$ is finite-dimensional, and since finite-dimensional C*-algebras are generated by a finite number of projections, it follows that this family equals $K(\CC(A))\cap\down C$. 
  Now let $F_1,F_2\in\F$. By Lemma~\ref{lem:joinindownC}, $C^*(F_1)\vee C^*(F_2)=C^*(F_1\cup F_2)$, making the family directed. Then:
  \[
    C 
    = \bigvee_{p\in P} C^*(\{p\}) = \mathop{\bigvee}\limits_{F\in\F} \bigvee_{p\in F} C^*(p) 
    = \mathop{\bigvee}\limits_{F\in\F} C^*(F),
  \]
  where the third equality used Lemma~\ref{lem:joinindownC}. Hence $\CC(A)$ is algebraic.
\end{proof}

\section{Quasi-continuity and quasi-algebraicity}\label{sec:quasi}

In this section we show that for dcpos $\CC(A)$ of C*-algebras $A$, the notions of quasi-continuity and quasi-algebraicity, which are generally weaker than continuity and algebraicity, are in fact equally strong. 

To define quasi-continuity and quasi-algebraicity \cite[Section III-3]{gierzetal:domains} we generalize the way below relation of a dcpo $\CC$ to nonempty subsets: write $\G \leq \HH$ when $\up \HH \subseteq \up \G$.
This is a pre-order, and we can talk about \emph{directed} families of nonempty subsets.
A nonempty subset $\G$ is \emph{way below} another one $\HH$, written $\G \ll \HH$, when $\bigvee \D \in \up \HH$ implies $D \in \up \G$ for some $D \in \D$.
Observe that $\{B\} \ll \{C\}$ if and only if $B \ll C$, so we may abbreviate $\G \ll \{C\}$ to $\G \ll C$, and $\{C\} \ll \HH$ to $C \ll \HH$. 

\begin{definition}\label{def:quasi}
  For an element $C$ in a dcpo $\CC$, define
  \begin{align*}
    \fin(C) & = \{ \F \subseteq \CC \mid \F \text{ is finite, nonempty, and } \F \ll C \}, \\
    \compfin(C) & = \{ \F\in\fin(C)\mid \F \ll \F \}.
  \end{align*}
  The dcpo is \emph{quasi-continuous} if each $\fin(C)$ is directed, and $C \nleq D$ implies $D \not\in \up \F$ for some $\F \in \fin(C)$. 
  It is \emph{quasi-algebraic} if each $\compfin(C)$ is directed, and $C \nleq D$ implies $D \not\in \up \F$ for some $\F \in \compfin(C)$.
\end{definition}

It is asserted without proof in \cite[Proposition III-3.10]{gierzetal:domains} that continuity implies quasi-continuity. For completeness let us sketch a proof. Let $\CC$ be continuous and $C\in\CC$. First one has to show that $\fin(C)$ is directed, so let $\F,\G\in\fin(C)$. Since $\ddown C$ is directed and has supremum $C$, it follows from $\F\ll C$ and $\G\ll C$ that there are $D_\F$ and $D_\G$ in $\ddown C$ such that $D_\F\in\up\F$ and $D_\G\in\up\G$. Since $\ddown C$ is directed, there is some $D\in\ddown C$ such that $D_\F,D_\G\leq D$. It now follows that $\mathcal H=\{D\}$ is an element in $\fin(C)$ such that $\F,\G\leq\mathcal H$. It follows by contraposition that $C\nleq D$ implies that $D\not\in\up\F$ for some $\F\in\fin(C)$. Indeed, if $D\in\up\F$ for each $\F\in\fin(C)$, then $D\in\up E$ for each $E\ll C$, which translates to $E\leq D$ for each $E\in\ddown C$. Hence $C=\bigvee\ddown C\leq D$ using the fact that $\CC$ is continuous. 
	
In almost the same way (replacing $\fin(C)$ by $\compfin(C)$, and replacing $\ddown C$ by the directed set of compact elements below $C$) one can show that algebraicity implies quasi-algebraicity.

Intuitively, quasi-continuity and quasi-algebraicity relax continuity and algebraicity to allow the information in approximants to be spread out over finitely many observations rather than be concentrated in a single one.

We start by analysing the way-below relation generalized to finite subsets.

\begin{lemma}\label{lem:fininCA}
  Let $A$ be a C*-algebra, $C\in\CC(A)$ and $\F\subseteq\CC(A)$. Then $\F\in\fin(C)$ if and only if $\F$ contains finitely many elements and $F\ll C$ for some $F\in\F$.
\end{lemma}
\begin{proof}
  Let $\F$ contain finitely many elements and assume that $F\ll C$ for some $F\in\F$. 
  Let $\D$ be a directed subset of $\CC(A)$ such that $C\subseteq\bigvee\D$. 
  Since $F\ll C$, we have $F\subseteq D$ for some $D\in\D$, so $D\in\up\F$. 
  Thus $\F\in\fin(C)$. 
  
  Conversely, suppose $\F\in\fin(C)$. Then $\F\ll C$ and $\F$ is nonempty and finite. 
  Now $\{C\}$ is a directed subset whose supremum contains $C$, so there is some $F = \{F_1,\ldots,F_n\} \in\F$ contained in $C$.
  Assume for a contradiction that each $F_i$ has infinite dimension. 
  Write $X$ for the spectrum of $C$, so $C\simeq C(X)$.
  Lemma~\ref{lem:topologicalconditionoffinitedimensionalsubalgebra} guarantees the existence of points $p_1,\ldots,p_n\in X$ with $F_j\nsubseteq C_{\overline{U_j}}$ for each open neighbourhood $U_j \subseteq X$ of $p_j$. 
  In particular, $F_j\nsubseteq\bigcap_{i=1}^nC_{\overline {U_i}}$ for each $i=1,\ldots,n$ and open neighbourhood $U_i$ of $p_i$. 
  Consider the family
  \[
    \Big\{ \bigcap_{i=1}^n C_{\overline{U_i}} \;\Big|\; U_i \text{ open neighbourhood of }p_i,\; i=1,\ldots,n\ \Big\}.
  \]
  It is directed and has supremum $C$ by Lemma~\ref{lem:directedfamilygeneratedbylocalneighbourhoods}.
  However, no member of the family contains $F_i$.
  If $F\in\F$ such that $F\nsubseteq C$, we cannot have $F\subseteq \bigcap C_{\overline{U_i}}$ for any $i$ or neighbourhood $U_i$ of $p_i$, because the latter is contained in $C$ by construction, contradicting $\F\ll C$. 
  We conclude that there must be a finite-dimensional $F\in\F$ such that $F\subseteq C$. Now $F \ll C$ follows from Proposition~\ref{prop:waybelowinCA}.
\end{proof}

\begin{lemma}\label{lem:compfinC}
  Let $A$ be a C*-algebra and let $C\in\CC(A)$.
  If $F \ll  C$, then $\{F\}\in\compfin(C)$. 
  If $\F\in\fin(C)$, then $\F \leq \F'$ for some $\F'\in\compfin(C)$.
\end{lemma}
\begin{proof}
  Let $F\ll C$. 
  By Lemma~\ref{lem:fininCA}, we have $\{F\}\in\fin(C)$. 
  By Lemma~\ref{prop:waybelowinCA}, we have $F\ll F$. 
  Therefore $\{F\}\ll\{F\}$, and so $\{F\}\in\compfin(C)$. 
  
  Let $\F\in\fin(C)$. 
  By Lemma~\ref{lem:fininCA}, there is an $F\in\F$ such that $F\ll C$. 
  The reasoning in the previous paragraph shows $\{F\}\in\compfin{C}$. 
  Since $F\in\F$, we have $F\in\up\F$, and so $\up\{F\}\subseteq\up\F$. We conclude that $\F\leq\F'$ for $\F'=\{F\}$.
\end{proof}

We are now ready to characterize quasi-continuity and quasi-algebraicity of $\CC(A)$.

\begin{theorem}
  A C*-algebra $A$ is scattered if and only if $\CC(A)$ is quasi-continuous, if and only if $\CC(A)$ is quasi-algebraic.
\end{theorem}
\begin{proof}
  If $A$ is scattered, then $\CC(A)$ is algebraic by Theorem~\ref{thm:algebraic}, and hence quasi-algebraic by the remarks following Definition~\ref{def:quasi}.

  Now assume that $\CC(A)$ is quasi-algebraic and let $C\in\CC(A)$. 
  Let $\F_1,\F_2\in\fin(C)$. 
  By Lemma~\ref{lem:compfinC}, there exist elements $\F_1',\F_2'\in\compfin(C)$ such that $\F_i\leq\F_i'$. 
  By quasi-algebraicity, $\compfin(C)$ is directed, so there is an $\F\in\compfin(C)$ such that $\F_1',\F_2'\leq\F$. 
  Hence $\F_1,\F_2\leq\F$. 
  Since $\compfin(C)\subseteq\fin(C)$, it then follows that $\fin(C)$ is directed. 
  Let $B\in\CC(A)$ satisfy $C\not\subseteq B$. Assume that $B\in\up\F$ for $\F\in\fin(C)$. 
  Lemma~\ref{lem:compfinC} provides $\F'\in\compfin(C)$ with $\F\leq\F'$. 
  But this means that $\up\F\subseteq\up\F'$. 
  Hence $B\in\up\F'$, which contradicts quasi-algebraicity. 
  Therefore we must have $B\notin\up\F$ for each $\F\in\fin(C)$, making $\CC(A)$ quasi-continuous.

  Finally, assume $\CC(A)$ is quasi-continuous. 
  Let $F_1,F_2\in\ddown C$. 
  By Lemma~\ref{lem:fininCA}, we have $\{F_1\},\{F_2\}\in\fin(C)$, and since $\fin(C)$ is directed, there is an $\F\in\fin(C)$ such that $\F\subseteq \up\{F_1\}\cap\up\{F_2\}$. 
  In other words, $F_1,F_2\subseteq F$ for each $F\in\F$, and since $\F\in\fin(C)$, Lemma~\ref{lem:fininCA} guarantees the existence of some $F$ such that $F\ll C$, making $\ddown C$ directed. 
  Let $B=\bigvee\ddown C$. Since $F\subseteq C$ for each $F\in\ddown C$, we have $B\subseteq C$. 
  If $B\neq C$, then $C\nsubseteq B$, so by quasi-continuity there must be an $\F\in\fin(C)$ with $B\notin\up\F$. 
  Hence $F\nsubseteq B$ for each $F\in\F$, and in particular, Lemma~\ref{lem:fininCA} implies the existence of some $F\in\F$ satisfying $F\ll C$, but $F\nsubseteq B$. 
  By definition of $B$ we have $F\subseteq B$ for each $F\ll C$, giving a contradiction. 
  Thus $\CC(A)$ is continuous, and by Theorem~\ref{thm:continuous}, $A$ is scattered.
\end{proof}

\section{Other notions of scatteredness}\label{sec:scattered}\label{sec:orderscattered}

In the previous sections, we have seen that the dcpo $\CC(A)$ of a C*-algebra $A$ is nice -- in the sense of being (quasi-)algebraic, (quasi-)continuous, or atomistic, which are all equivalent -- precisely when $A$ is scattered. 
In this section we study when $\CC(A)$ itself is scattered, in two different ways: by putting a topology on $\CC(A)$ and asking when it is a scattered topological space; and by considering an established (but different) notion of scatteredness for partially ordered sets directly on $\CC(A)$. Both will turn out to be very restrictive, in the sense that they coincide with $A$ being finite-dimensional.

If the C*-algebra $A$ is scattered, then we can turn the domain $\CC(A)$ itself into the spectrum of another C*-algebra, which this section studies. We will use the Lawson topology, that turns approximation in domains into topological convergence.

\begin{definition}\label{def:scottlawson}
  The \emph{Scott topology} declares subsets $\U$ of a dcpo to be open if $\up \U = \U$, and $\D \cap \U \neq \emptyset$ when $\bigvee^{\up} \D \in \U$.
  The \emph{Lawson topology} has as basic open subsets $\U \setminus \up \F$ for a Scott open subset $\U$ and a finite subset $\F$.
\end{definition}

These topologies capture approximation in the following sense. A monotone function $f$ between dcpos is continuous with respect to the Scott topology precisely when it is Scott continuous, i.e., when $\bigvee f[\D]=f(\bigvee\D)$ for directed subsets $\D$~\cite[Proposition II-2.1]{gierzetal:domains}. 
Thus Proposition~\ref{prop:CfisScottcontinuous} shows that $\CC(f)$ is Scott continuous.
Similarly, a meet-semilattice homomorphism between complete semilattices is continuous with respect to the Lawson topology precisely when $\bigwedge f[\D]=f(\bigwedge\D)$ for nonempty subsets $\D$~\cite[Theorem III-1.8]{gierzetal:domains}.

\begin{proposition}
  If $A$ is a scattered C*-algebra, then $\CC(A)$ is a totally disconnected compact Hausdorff space in the Lawson topology.
\end{proposition}
\begin{proof}
  If $A$ is scattered, then $\CC(A)$ is both an algebraic domain and a complete semilattice.
  Therefore it is compact Hausdorff in the Lawson topology~\cite[Corollary III-1.11]{gierzetal:domains}. 
  Moreover, it follows that $\CC(A)$ is zero-dimensional~\cite[Exercise III-1.14]{gierzetal:domains}, which for compact Hausdorff spaces is equivalent to being totally disconnected \cite[Theorem 29.7]{willard:generaltopology}.
\end{proof}

For example, let $A$ be the algebra of all $2$-by-$2$-matrices as discussed in Example \ref{ex:BH}. In $\CC(A)$ there is then one bottom element, and all other elements are incomparable to each other and are in 1-1 correspondence with bases of $\C^2$. As a consequence, the Lawson topology of $\CC(A)$ is homeomorphic to the one-point compactification of a discrete space of cardinality $2^{\aleph_0}$.

It follows from the previous proposition that any scattered C*-algebra $A$ gives rise to another, commutative, C*-algebra $C(X)$ for $X=\CC(A)$ with its Lawson topology. 
Thus we can speak about the domain of commutative C*-subalgebras entirely within the language of C*-algebras.
However, there is a caveat:
iterating this construction only makes sense in the finite-dimensional case, as the following theorem shows.

\begin{theorem}
  A scattered C*-algebra $A$ is finite-dimensional if and only if $\CC(A)$ is scattered in the Lawson topology.
\end{theorem}
\begin{proof}
  Let $A$ have finite dimension, so it is certainly scattered, and $\CC(A)$ is algebraic. 
  It follows that a basis for the Scott topology is given by $\up C$ for $C$ compact~\cite[Corollary II-1.15]{gierzetal:domains}. 
  Thus sets of the form $\up C\setminus\up\F$ with $C$ compact and $\F$ finite form a basis for the Lawson topology. 
  Take a nonempty subset $\mathcal{S} \subseteq \CC(A)$, and let $M$ be a maximal element of $\mathcal{S}$, which exists by~\cite[Lemma 10]{lindenhovius:artinian}.
  Since $M$ must be finite-dimensional too, it is compact by Proposition~\ref{prop:finitedimiscompact}.
  Hence $\up M$ is Scott open and therefore Lawson open. 
  Maximality of $M$ in $\mathcal{S}$ now gives $\mathcal{S} \cap\up M=\{M\}$, and since $\up M$ is Lawson open, it follows that $M$ is an isolated point of $\mathcal S$. Hence $\CC(A)$ is scattered (cf. Definition \ref{def:topologically scattered}). 
  
  For the converse, assume $A$ is infinite-dimensional. 
  Then $\CC(A)$ has a noncompact element $C$, and $\down C$ contains an isolated point if it intersects some basic Lawson open set in a single point.
  Hence $\down C \,\cap\, \up K \,\setminus\, \up\F$ must be a singleton for some finite set $\F\subseteq\CC(A)$ and some compact $K \in \CC(A)$. 
  In other words, $[K,C]\setminus\up\F$ is a singleton, where $[K,C]$ is the interval $\{D\in\CC(A) \mid K\subseteq D\subseteq C\}$. 
  Since $C$ is infinite-dimensional and scattered (by Theorem~\ref{thm:Kusuda}), there are infinitely many atoms in $[K,C]$: for $\CC(A)$ is atomistic by Theorem~\ref{thm:CAatomistic} and hence $C$ dominates infinitely many atoms $C_i$, but $K$ is finite-dimensional by Proposition~\ref{prop:finitedimiscompact}, so that $C_i \vee K$, excepting the finitely many $C_i \leq K$, give infinitely many atoms in $[K,C]$.
  Hence there is no finite subset $\F\subseteq\CC(A)$ making $[K,C]\setminus\up\F$ a singleton. 
  We conclude that $\down C$ has no isolated points, so $\CC(A)$ cannot be scattered. 
\end{proof}

Just like there is an established notion of scatteredness for topological spaces and C*-algebras, there is an established notion of scatteredness for partially ordered sets. In the rest of this section we show that the two notions diverge, and should not be confused.

\begin{definition}
  A chain $C$ in a poset $P$ is \emph{order-dense} if none of its elements covers another one, \ie if $x<z$ in $C$ then $x<y<z$ for some $y\in C$.
  A poset is \emph{order-scattered} when it does not contain an order-dense chain of at least two points.
\end{definition}

\begin{lemma}\label{lem:notorderscattered}
  If a C*-algebra $A$ is not scattered, then $\CC(A)$ is not order-scattered.
\end{lemma}
\begin{proof}
  If $A$ is not scattered it has a commutative C*-subalgebra with spectrum $[0,1]$ by Theorem~\ref{thm:Kusuda}, so without loss of generality we may assume $A=C([0,1])$. Consider 
  \[
    \{ C_{[x,1]} \mid x \in [0,1) \} \subseteq \CC(A).
  \]
  Because $C_{[x,1]} \subseteq C_{[y,1]}$ if and only if $x \leq y$, this set is order-isomorphic to the order-dense chain $[0,1)$ via the map $C_{[x,1]} \mapsto x$.
\end{proof}

\begin{lemma}\label{lem:orderscattered}
  Let $X$ be an infinite scattered compact Hausdorff space, and let $A=C(X)$. Then $\CC(A)$ is not order-scattered.
\end{lemma}
\begin{proof}
  First observe that $X$ must contain an infinite number of isolated points, for if $X$ had only finitely many isolated points $x_1,\ldots,x_n$, then $X \setminus \{x_1,\ldots,x_n\}$ would be closed and hence contain an isolated point $x_{n+1}$, which is also isolated in $X$ because $X \setminus \{x_1,\ldots,x_n\}$ would be open.
  Choose a countably infinite set $Y$ of isolated points of $X$. Note that $Y$ is open, but cannot be closed because $X$ is compact.
  Let $Z = \overline{Y} \setminus Y$ be the boundary of $Y$. Since $Y$ is open, $Z$ is closed. Moreover, if $S \subseteq Y$, then $Z \cup S = (\overline{Y} \setminus Y) \cup S = \overline{S} \setminus (Y \setminus S)$ is closed because $Y \setminus S$ consists only of isolated points and hence is open. As $Y$ is countably infinite, we can label its elements by rational numbers $Y = \{x_q\}_{q \in \mathbb{Q}}$. For each $q \in \mathbb{Q}$, set
  \[
    K_q = Z \cup \{x_r \mid r \leq q\}\text{,}
  \]
  and notice that $K_q$ is closed and infinite. Now $q \mapsto K_q$ is an order embedding of $\mathbb{Q}$ into the set $\mathcal{F}(X)$ of all closed subsets of $X$ with at least two points, partially ordered by inclusion. In turn, $K \mapsto C_K$ is an order embedding of $\mathcal{F}(X)\op$ into $\CC(A)$. Composing gives an order embedding $\mathbb{Q}\op \to \CC(A)$, and therefore an order-dense chain in $\CC(A)$.
\end{proof}

\begin{theorem}\label{thm:orderscattered}
  A C*-algebra $A$ is finite-dimensional if and only if $\CC(A)$ is order-scattered.
\end{theorem}
\begin{proof}
  If $A$ is finite-dimensional, then so is each $C \in \CC(A)$. Hence all chains in $\CC(A)$ have finite length, and therefore cannot be order-dense. That is, $\CC(A)$ is order-scattered.

  For the converse, assume that $A$ is infinite-dimensional. 
  We distinguish two cases. If $A$ is not scattered, then Lemma~\ref{lem:notorderscattered} shows that $\CC(A)$ is not order-scattered.
  If $A$ is scattered, then it has a maximal commutative C*-subalgebra with scattered spectrum $X$. Because $C(X)$ must be infinite-dimensional~\cite[Exercise 4.6.12]{kadisonringrose:oa3}, $X$ is infinite. Now Lemma~\ref{lem:orderscattered} shows that $\CC(C(X))$, and hence $\CC(A)$, is not order-scattered.
\end{proof}

\section{Projections and posets of Boolean subalgebras}\label{sec:Booleansubalgebras}

We recall that an element $p$ of a C*-algebra $A$ is called a projection when $p^2=p=p^*$. The set $\Proj(A)$ of all projections in $A$ can be ordered where $p \leq q$ if and only if $p=pq$. In general, if $p$ is a projection then so is $1-p$, and in fact the projections form an orthomodular poset, which we define below. 

In this section we aim to reconstruct the orthomodular poset $\Proj(A)$ of projections in a C*-algebra $A$ from $\CC(A)$ for the reason that in many cases $\Proj(A)$ encodes much of the structure of $A$, especially if $A$ belongs to a class of C*-algebras that have an ample supply of projections such as AF-algebras or von Neumann algebras (cf.\ Definition \ref{def:vN algebra} below). Given an orthomodular poset $P$, one can consider the poset $\B(P)$ of Boolean subalgebras of $P$, which we will also define below, and which already had been proven to determine $P$ up to isomorphism~\cite{hardingheunenlindenhoviusnavara:booleansubalgebras}. We aim to exploit this fact, hence we will introduce the following subposet of $\CC(A)$ that turns out to be isomorphic to $\B(\Proj(A))$:

\begin{definition}
	Let $A$ be a C*-algebra. Then we denote the subposet of $\CC(A)$ consisting of the commutative C*-subalgebras of $A$ that are AF-algebras by $\CC_\AF(A)$.
\end{definition}

\begin{lemma}
	Let $A$ be a C*-algebra. Then $\CC_\AF(A)$ is a dcpo.
\end{lemma}	 
\begin{proof}
	Let $\D\subseteq\CC_\AF(A)$ be directed, and consider its supremum $S=\overline{\bigcup\D}$ in $\CC(A)$. We will show that $S\in\CC_\AF(A)$. Let $a_1,\ldots,a_n\in\bigvee\D$ and let $\varepsilon>0$. Then there are $d_1,\ldots,d_n\in\bigcup\D$ such that $\|a_i-d_i\|<\varepsilon/2$. Let $D_i\in\D$ such that $d_i\in D_i$. Since $\D$ is directed, there is a $D$ in $\D$ containing $D_1,\ldots,D_n$, hence containing $d_1,\ldots,d_n$. Since $D\in\CC_\AF(A)$, Proposition \ref{prop:commutativeAFistotallydisconnected} assures that there is a finite-dimensional C*-subalgebra $B\subseteq D$ and $b_1,\ldots,b_n\in B$ such that $\|d_i-b_i\|<\varepsilon/2$. Hence \[\|a_i-b_i\|<\|a_i-d_i\|+\|d_i-b_i\|<\varepsilon,\] hence Proposition \ref{prop:commutativeAFistotallydisconnected} assures that $S$ is an AF-algebra.
\end{proof}
It follows from the next lemma that $\CC_\AF$ is a functor $\CStar\to\DCPO$:
\begin{lemma}
	Let $f:A\to B$ be a $*$-homomorphism between C*-algebras $A$ and $B$. Then $\CC(f):\CC(A)\to\CC(B)$ restricts to a Scott continuous map $\CC_\AF(A)$ and $\CC_\AF(B)$. 
\end{lemma}
\begin{proof}
  First consider the case $B=0$. Then $B$ is the terminal object of the category of unital C*-algebras, so $f$ is the unique $*$-homomorphism $A\to B$. Now $\CC(f)$ is a monotone map to the $1$-element poset $\CC(B)$, and must therefore be Scott continuous. Hence we may assume $B\neq 0$.	
  Since $f$ is a $*$-homomorphism, it is bounded, and it is non-zero since $f(1)=1$, so $\|f\|\neq 0$. Let $C\in\CC_\AF(A)$. Then $\CC(f)(C)=f[C]$, hence let $x_1,\ldots,x_n\in f[C]$, i.e., there are $c_1,\ldots,c_n\in C$ such that $x_i=f(c_i)$ for each $i\in\{1,\ldots,n\}$. Let $\varepsilon>0$. Since $C$ is approximately finite-dimensional and commutative, it follows from Proposition \ref{prop:commutativeAFistotallydisconnected} that there is a finite-dimensional C*-algebra $D\subseteq C$ containing $d_1,\ldots, d_n$ such that $\|c_i-d_i\|<\varepsilon/\|f\|$ for each $i\in\{1,\ldots,n\}$. Note that since $f$ is linear, $f[D]$ must be a finite-dimensional C*-subalgebra of $f[C]$. Hence 
	\[\|f(c_i)-f(d_i)\|\leq\|f\|\|c_i-d_i\|<\varepsilon,\] which, in combination with Proposition \ref{prop:commutativeAFistotallydisconnected}, shows that $f[C]$ is an AF-algebra.
	Hence $\CC(f)$ restricts to a map between $\CC_\AF(A)$ and $\CC_\AF(B)$, and since $\CC(f)$ is Scott continuous (cf. Proposition \ref{prop:CfisScottcontinuous}), so is its restriction.
\end{proof}

Next we define orthomodular posets. For a more detailed overview of orthomodular structures we refer to \cite{dvurecenskijpulmannova:quantumstructures}.
\begin{definition}
  A partially ordered set $P$ is an \emph{orthoposet} when it has a greatest element $1$ and a least element $0$, and it comes equipped with an operation $\perp \colon P \to P$, called the \emph{orthocomplementation}, satisfying for each $p,q \in P$:
		\begin{itemize}
			\item $p^{\perp\perp}=p$;
			\item if $p \leq q$, then $q^\perp \leq p^\perp$;
			\item $p$ and $p^\perp$ have a least upper bound and greatest upper bound, and $p\wedge p^\perp=0$ and $p\vee p^\perp=1$.
		\end{itemize}
		If $p\leq q^\perp$ (or equivalently $q\leq p^\perp$), then we say that $p$ and $q$ are \emph{orthogonal} and write $p\perp q$.
		If $P$ is an orthoposet for which 
		\begin{itemize}
			\item $p \perp q$ implies the existence of $p\vee q$;
			\item $p\leq q$ implies $q=p\vee(q\wedge p^\perp)$,
		\end{itemize}
		then we call $P$ an \emph{orthomodular poset}. If, in addition, $P$ is a lattice, we call it an \emph{orthomodular lattice}.
			If $p$ and $q$ are elements in an orthomodular poset $P$ for which there are pairwise orthogonal elements $e_1,e_2,e_3\in P$ such that
			\begin{align*}
			p&=e_1\vee e_3, &  q& =e_2\vee e_3,
			\end{align*}
			then we say that $p$ and $q$ \emph{commute}. 
\end{definition}
If $A$ is a C*-algebra, then $\Proj(A)$ becomes an orthomodular poset if we define its orthocomplementation by $p^\perp=1-p$. Note that there are C*-algebras $A$ for which $\Proj(A)$ is not a lattice, see for instance \cite[Lemma 2.1]{lazar:afalgebraslattice}. It is easy to see two projections $p$ and $q$ in a C*-algebra $A$ are orthogonal in the orthomodular poset $\Proj(A)$ if and only if they are orthogonal in the operator-algebraic sense $pq=0$. Similarly, $p$ and $q$ commute in the orthomodular poset $\Proj(A)$ if and only if they commute in an algebraic sense: $pq=qp$.

Next we define the appropriate morphisms for orthomodular posets:
\begin{definition}
		
			A map $f:P\to Q$ between orthomodular posets $P$ and $Q$ is called an \emph{orthomodular morphism} if 
			\begin{itemize}
				\item $f(1)=1$;
				\item $f(x^\perp)=f(x)^\perp$ for each $x\in P$;
				\item $x\perp y$ implies $f(x\vee y)=f(x)\vee f(y)$ for each $x,y\in P$.
			\end{itemize}
			We denote the category of orthomodular posets with orthomodular morphisms by $\OMP$.
\end{definition}
Let $f:A\to B$ be a $*$-homomorphism between C*-algebras $A$ and $B$. Since it preserves all algebraic operations, it follows that $f(p)$ is a projection in $B$ if $p$ is a projection in $A$, hence $f$ restricts to a map $\Proj(A)\to\Proj(B)$. If we define $\Proj(f)$ to be the restriction of $f$ to $\Proj(A)$, it is routine to check that $\Proj$ becomes a functor $\CStar\to\OMP$.

If $P$ is a Boolean algebra, then any pair $x,y$ of elements in $P$ commute, since we can write $x=e_1\vee e_3$ and $y=e_2\vee e_3$ with $e_1=x\wedge y^\perp$, $e_2=x^\perp\wedge y$ and $e_3=x\wedge y$. Conversely, if $P$ is an orthomodular poset for which all elements mutually commute, then $P$ is a Boolean algebra, which follows from the Foulis--Holland Theorem \cite[Theorem 12.3.1]{beltrametticassinelli:logicquantum}. Hence we can regard Boolean algebras as `commutative' orthomodular posets, which gives more reason to consider the poset $\B(P)$ of Boolean subalgebras of an orthomodular poset, which we define now in more detail.

\begin{definition}\label{def:B}
  Let $P$ be an orthomodular poset. A subset $B$ that is closed under the operation $x\mapsto x^\perp$ and for which the join of any two mutually orthogonal elements is contained in $B$ is called a \emph{sub-orthomodular poset}, which becomes an orthomodular poset if we equip it with the order and the orthocomplementation inherited from $P$. If, in addition, $B$ is a Boolean algebra, i.e., it is an orthocomplemented lattice in which the distributive law holds: for each $x,y,z\in B$, we have
  \[x\wedge(y\vee z)=(x\wedge y)\vee(x\wedge z),\]
   then we call $B$ a \emph{Boolean subalgebra} of $P$. We denote the set of all Boolean subalgebras of $P$ by $\B(P)$, which we partially order by inclusion.
\end{definition}

\begin{proposition}\label{prop:BofPisalgebraiccompletesemilattice}
	For each orthomodular poset $P$, the poset $\B(P)$ is an algebraic semilattice, where $\bigvee\D=\bigcup\D$ for each directed set $\D\subseteq\B(P)$, and $\bigwedge\mathcal S=\bigcap\mathcal S$ for any non-empty subset $\mathcal S\subseteq\B(P)$. Moreover, the compact elements of $\B(P)$ are precisely the finite Boolean subalgebras.
\end{proposition}	
\begin{proof}
	Let $\mathcal S$ be a non-empty collection of Boolean subalgebras of $P$, and let $B$ be its intersection. Let $x,y\in B$ (we do not need the assumption that $x$ and $y$ are orthogonal). Then $x^\perp,x\vee y\in D$ for each $D\in\mathcal S$, so $x^\perp,x\vee y\in B$, which shows that $B$ is a sub-orthomodular poset of $P$. Since the distributivity law holds in each $D\in\mathcal S$, it follows that it holds in $B$, hence $B$ is a Boolean subalgebra of $P$. Clearly $B$ is the infimum of $\mathcal S$.
	
	Now let $\D$ be a directed set in $\B(P)$ and let $B$ be its union. Let $x,y\in B$. Then there are $D_x,D_y\in\D$ such that $x\in D_x$ and $y\in D_y$, and since $\D$ is directed there is some $D\in\D$ such that $D_x,D_y\subseteq D$. Hence $x,y\in D$, which is a Boolean subalgebra, hence $x\vee y$ and $x^\bot$ exist and are contained in $D$, whence they are elements of $B$. In a similar way, if $x,y,z\in B$, we can find a $D\in\D$ containing $x$, $y$ and $z$, and since $D$ is Boolean, $x$, $y$, $z$ satisfy the distributivity law. We conclude that $B$ is a Boolean subalgebra of $P$. Clearly $B$ is the supremum of $\D$.
	
	Let $B$ be a finite Boolean subalgebra, and write $B=\{b_1,\ldots,b_n\}$. Let $\D\subseteq\B(P)$ directed such that $B\subseteq\bigvee\D$. Since $\bigvee\D=\bigcup\D$, it follows that for each $i\in\{1,\ldots,n\}$ there is a $D_i\in\D$ such that $b_i\in D_i$. Since $\D$ is directed, there is some $D\in\D$ such that $D_1,\ldots,D_n\subseteq D$. Hence $B\subseteq D$, which shows that $B$ is compact.
	
	For the converse, we first assume that $B$ is an arbitrary Boolean subalgebra of $P$, and let $\D$ be the set of all finite Boolean subalgebras of $B$. Then $\D$ is directed: if $D_1,D_2\in\D$, let $S=D_1\cup D_2$. Since $S$ is finite, so is $S\cup S^\perp$, where $S=\{s^\perp:s\in S\}$. Then consider $R=\{\bigwedge F:F\subseteq S\cup S^\perp\}$, which is also finite since $S\cup S^\perp$ is finite. Finally, let $D=\{\bigvee X:X\subseteq R\}$ which is finite for $R$ is finite. It now follows from the De Morgan Laws that $D$ is a subset of $B$ that is closed under meets, joins and the orthocomplementation, whence $D$ is a Boolean subalgebra. For each $b\in B$, the set $\{0,b,b^\perp,1\}$ forms a finite Boolean subalgebra, hence it follows that $B=\bigcup\D$, i.e., $B=\bigvee\D$. This shows that all Boolean subalgebras of $P$ are the directed supremum of finite Boolean subalgebras. Moreover, if $B$ is compact, then it follows that $B\subseteq D$ for some finite Boolean subalgebra $D$ of $P$, whence $B$ is finite, too, which concludes the proof that the finite Boolean subalgebras of $P$ are exactly the compact elements of $\B(P)$. Thus each $B\in\B(P)$ is the directed supremum of compact elements, whence $\B(P)$ is algebraic.
\end{proof}	

The previous proposition generalizes the statement that the lattice of subalgebras a Boolean algebra is algebraic~\cite{graetzerkohmakkai:booleansubalgebras} from Boolean algebras to arbitrary orthomodular posets. See also~\cite{heunen:piecewiseboolean} for a different generalization. 

To show that the assignment $P\mapsto\B(P)$ extends to a functor $\OMP\to\DCPO$, we first need a lemma. 
\begin{lemma}\label{lem:ompmorphism}
	Let $f:P\to Q$ be an orthomodular morphism between orthomodular posets $P$ and $Q$, then $f$ preserves binary joins of of commuting elements. 
\end{lemma}
\begin{proof}
	Since $x$ and $y$ commute, we have $x=e_1\vee e_3$, $y=e_2\vee e_3$ and $x\vee y=e_1\vee e_2\vee e_3$ where $e_1,e_2,e_3$ are mutually orthogonal, whence
	\begin{align*}
	f(x\vee y) & =f(e_1\vee e_2\vee e_3)=f(e_1\vee e_2)\vee f(e_3)=f(e_1)\vee f(e_2)\vee f(e_3)\\
	& =f(e_1)\vee f(e_3)\vee f(e_2)\vee f(e_3)=
	f(e_1\vee e_3)\vee f(e_2\vee e_3)=f(x)\vee f(y). 
	\end{align*} 
  Thus $f$ preserves binary joins of commuting elements.
\end{proof}

\begin{proposition}\label{prop:Bisfunctor}
	Let $P$ and $Q$ be orthomodular posets and let $f:P\to Q$ be an orthomodular morphism. Then $B \mapsto f[B]$ is a well-defined map $\B(f) \colon \B(P)\to\B(Q)$ that is Scott continuous. Moreover, if $f$ is injective, then $\B(f)$ is an order embedding.
\end{proposition}
\begin{proof}
	Let $B\in\B(P)$, and let $x,y\in B$. Since $B$ is Boolean, $x$ and $y$ commute. Lemma \ref{lem:ompmorphism} assures that $\varphi$ preserves their join: $f(x\vee y)=f(x)\vee f(y)$. Thus the join of $f(x)$ and $f(y)$ exists and is contained in $f[B]$. By definition of an orthomodular morphism, $f$ preserves orthocomplementation, so $f(x)^\perp\in f[B]$. Since the join of all elements in $B$ exist, and De Morgan's laws hold in orthomodular posets:
	\[f(x\wedge y)=f((x^\perp\vee y^\perp)^\perp)=f(x^\perp\vee y^\perp)^\perp=(f(x)^\perp\vee f(y)^\perp)^\perp=f(x)\wedge f(y)\text.\] Therefore the meet of $f(x)$ and $f(y)$ exists and is contained in $f[B]$, and $f$ preserves all binary meets. Since $B$ is Boolean, it satisfies the distributive law, and since $f$ preserves all binary meets, binary joins and the orthocomplementation, it follows that $f[B]$ also satisfies the distributive law. We conclude that $f[B]$ is a Boolean subalgebra of $Q$, which shows that $\B(f)$ is well-defined. Clearly, it preserves inclusions, and since the direct image preserves unions, and the supremum of a directed family in $\B(P)$ is given by the union of its members (cf.\ Proposition~\ref{prop:BofPisalgebraiccompletesemilattice}), it follows that $\B(f)$ is Scott continuous.  

	Finally, assume that $f$ is injective, and let $B_1,B_2\in\B(P)$ such that $\B(f)(B_1)\subseteq\B(f)(B_2)$. This implies that $f[B_1]\subseteq f[B_2]$. Let $x\in B_1$. Then $f(x)\in f[B_2]$, hence there is some $y\in B_2$ such that $f(x)=f(y)$. By injectivity of $f$ it follows that $x=y$, whence $x\in B_2$. 
\end{proof}

\begin{theorem}\label{thm:CAFisBProjA}
	$\CC_\AF$ and $\B\circ\Proj$ are functors $\CStar\to\DCPO$ that are naturally isomorphic. In particular, if $A$ is a C*-algebra, then $\CC_\AF(A) \simeq \B(\Proj(A))$ via $C \mapsto \Proj(C)$.
\end{theorem} 
  \begin{proof}
    Let $A$ be a C*-algebra. Write $\varphi_A\colon \CC_\AF(A) \to \B(\Proj(A))$ for the map $C \mapsto \Proj(C)$.
    Define $\psi_A \colon \B(\Proj(A)) \to \CC_\AF(A)$ by $\psi_A(B) = C^*(B)$; the proof of Proposition~\ref{prop:elementsofCgeneratedbyprojections} shows that this is well-defined. 
    Both $\varphi_A$ and $\psi_A$ are clearly monotone. Moreover, if $C \in \CC_\AF(A)$, then $C^*(\Proj(C))=C$, so that $\psi_A(\varphi_A(C))=C$.
    Now let $B \in \B(\Proj(A))$. Say that $B$ is isomorphic to the Boolean algebra $B(X)$ of clopen subsets of the Stone space $X$. There is an  isomorphism $B(X)\simeq\Proj(C(X))$. Hence we may assume that $B=\Proj(C(X))$ for some Stone space $X$. 
    Now $C^*(B)=C^*(\Proj(C(X)))=C(X)$, whence $\Proj(C^*(B))=B$, so that $\varphi_A(\psi_A(B))=B$.
    Therefore $\varphi_A$ and $\psi_A$ are inverses, hence order isomorphisms.
    In order to show that the functors are naturally isomorphic, we show that $\varphi_A$ forms the $A$-component of a natural isomorphism $\CC_\AF\to\B\circ\Proj$. So let $f:A\to A'$ be a $*$-homomorphism between C*-algebras $A$ and $A'$. We show that the following diagram commutes:
    \begin{equation*}
    \xymatrix{\CC_\AF(A)\ar[dd]_{\CC_\AF(f)}\ar[rr]^{\varphi_A}  && \B(\Proj(A))\ar[dd]^{\B(\Proj(f))}  \\
    	\\
    	\CC_\AF(A')\ar[rr]_{\varphi_{A'}} && \B(\Proj(A')).}
    \end{equation*}
     Let $C\in\CC_\AF(A)$. Then it follows from Proposition \ref{prop:commutativeAFistotallydisconnected} that $C=C^*(\Proj(C))$, hence
     \begin{align*}
     \varphi_{A'}\circ\CC_\AF(f)(C) & = \varphi_{A'}(f[C]) = \varphi_{A'}(f[C^*(\Proj(C))])= \varphi_{A'}(C^*(f[\Proj(C)]))\\ & = \varphi_{A'}\circ\psi_{A'}(f[\Proj(C)]) = f[\Proj(C)] = \Proj(f)[\Proj(C)]\\
     &= \B(\Proj(f))(\Proj(C)) = \B(\Proj(f))\circ\varphi_A(C),
     \end{align*}
     where the first equality follows by definition of $\CC_\AF(f)$ as a restriction of $\CC(f)$, the third equality follows from Lemma \ref{lem:imageofgeneratedsubalgebra}, the fourth equality by definition of $\psi_{A'}$, which is the inverse of $\varphi_{A'}$, which gives the fifth equality, and the penultimate equality follows by definition action of $\Proj$ on $*$-homomorphisms.
\end{proof}

It follows that $\CC_\AF(A)$ is an algebraic complete semilattice, whose compact elements are precisely the finite-dimensional commutative C*-subalgebras~\cite[Corollary 6.2.5]{lindenhovius:thesis}. See also~\cite{heunen:piecewiseboolean}.
 
Combining the previous theorem with~\cite{hardingheunenlindenhoviusnavara:booleansubalgebras} gives the following.

\begin{corollary}\label{cor:reconstructionProj}
	We can construct an orthomodular poset $P$ orthomodular isomorphic to $\Proj(A)$ completely in terms of $\CC_\AF(A)$.
\end{corollary}

The next lemma characterizes $\CC_\AF(A)$ as a subposet of $\CC(A)$.

\begin{lemma}\label{lem:characterizationCCAF}
	Let $A$ be a C*-algebra. Then $\CC_\AF(A)$ is the subposet of $\CC(A)$ consisting of all elements that are the supremum of some subset of atoms of $\CC(A)$.
	\end{lemma}
	\begin{proof}
This follows directly from combining Proposition \ref{prop:commutativeAFistotallydisconnected} and Proposition \ref{prop:elementsofCgeneratedbyprojections}.
		\end{proof}
	
Combining Corollary \ref{cor:reconstructionProj} and Lemma \ref{lem:characterizationCCAF} now gives the main result of this section.

\begin{theorem}
  Let $A$ be a C*-algebra. Then we can construct an orthomodular poset $P$ isomorphic to $\Proj(A)$ completely in terms of $\CC(A)$.
\end{theorem}

\section{AW*-algebras}\label{sec:projections}

The previous section observed that the projections of a commutative C*-algebra form a Boolean algebra. In particular, the projections of $C(X)$ are precisely the indicator functions of clopen subsets of $X$. 
Thus, if a C*-algebra has many projections, it is intuitively rather disconnected. As an example, we have seen in Proposition \ref{prop:commutativeAFistotallydisconnected} that every commutative AF-algebra has a Gelfand spectrum that is totally disconnected, or a Stone space, since every Gelfand spectrum is also compact Hausdorff. Hence we can regard AF-algebras as operator-algebraic versions of Stone spaces. A stronger notion is a disconnected space is that of an \emph{extremally disconnected} space, in which the closure of an open set is open. Extremally disconnected compact Hausdorff spaces are also called \emph{Stonean} spaces. The operator algebraic version is as follows~\cite{kaplansky:projections}.

\begin{definition}\label{def:AW-algebra}
	An \emph{AW*-algebra} is a C*-algebra $A$ such that $\Proj(A)$ is a complete lattice, and every maximal element in $\CC(A)$ is generated by its projections.
\end{definition}

Equivalently, but more in line with our purposes, an AW*-algebra is a C*-algebra whose maximal C*-subalgebras have a Stonean Gelfand spectrum~\cite[Theorem~8.2.5]{saitowright:monotone}. The most prominent examples of AW*-algebras are the von Neumann algebras, which we define now.
\begin{definition}\label{def:vN algebra}
Let $H$ be a Hilbert space, and let $S\subseteq B(H)$. Then the \emph{commutant} $S'$ of $S$ is the set $\{b\in B(H):ab=ba\text{ for each }a\in S\}$. If $M\subseteq B(H)$ is a C*-algebra that it equal to its \emph{bicommutant} $M''$, then we call $M$ a \emph{von Neumann algebra} on $B(H)$. A C*-algebra *-isomorphic to some von Neumann algebra is called a \emph{W*-algebra}. 	
	\end{definition}
\begin{example}\cite[Theorem III.1.18]{takesaki:oa1}\label{ex:l infty}
Let $X$ be a set. Then $\ell^\infty(X)$, the space of all complex-valued functions $f$ on a set $X$ such that $\sup_{x\in X}|f(x)|$ is a finite number (which we define to be the norm of $f$), is a commutative W*-algebra.  	
\end{example}	

In this section we consider variations on $\CC(A)$ that cooperate well with projections. We start with approximating an AW*-algebra by its commutative AW*-subalgebras. 

\begin{definition}\label{def:AofA}
	An AW*-subalgebra of an AW*-algebra $A$ is a C*-subalgebra $C \subseteq A$ that is an AW*-algebra in its own right, with the same suprema of projections as in $A$.  Write $\A(A)$ for the partially ordered set of commutative AW*-subalgebras of $A$ under inclusion.
\end{definition}

There is a similar notion $\V(A)$ of W*-subalgebras of a W*-algebra $A$, that is studied in~\cite{doeringbarbosa:domain}. Here a W*-subalgebra of the W*-algebra $A$ is just an AW*-subalgebra, hence $\V(A)=\A(A)$. Thus working in the setting of AW*-algebras is more general. We note that if $M\subseteq B(H)$ is a von Neumann algebra, then a C*-subalgebra $N$ of $M$ is a W*-subalgebra if and only if it is a von Neumann algebra on $B(H)$ (cf. \cite[Exercise 4.24]{berberian}). Hence $\V(M)=\{C\subseteq M \mid C\text{ is a commutative von Neumann algebra on }B(H)\}.$

\begin{proposition}\label{prop:galois}
	For an AW*-algebra $A$ there is a Galois correspondence
	\[\begin{tikzpicture}[xscale=1.5]
	\node (v) at (-1,0) {$\A(A)$};
	\node (c) at (1,0) {$\CC(A)$};
	\node at (0,0) {$\perp$};
	\draw[->] ([yshift=2mm]c.west) to ([yshift=2mm]v.east);
	\draw[right hook->] ([yshift=-2mm]v.east) to ([yshift=-2mm]c.west);
	\end{tikzpicture}\]  
	where the upper adjoint maps a C*-subalgebra $C \in \CC(A)$ to the smallest AW*-subalgebra of $A$ containing it.
	Hence $\A(A)$ is a dcpo.
\end{proposition}
\begin{proof}
	Write $C' = \bigcap \{W \in \A(A) \mid C \subseteq W\}$ for the smallest AW*-subalgebra of $A$ containing $C \in \CC(A)$, which exists by \cite[Proposition 4.8.(ii)]{berberian}.
	If $C \subseteq D$, then clearly $C' \subseteq D'$.
	By construction we have $C' \subseteq W$ if and only if $C \subseteq W$, for $C \in \CC(A)$ and $W \in \A(A)$.
	Finally, notice that $W'=W$ for $W \in \A(A)$.
\end{proof}

Next we characterize the AW*-algebras $A$ whose dcpos $\CC(A)$ and $\A(A)$ are continuous, extending~\cite[Theorem 6.1]{doeringbarbosa:domain}.
This needs the following lemma.

\begin{lemma}\label{lem:stonescattered}
	Compact Hausdorff spaces that are scattered and Stonean must be finite. 
\end{lemma}
\begin{proof}
	Consider the open and discrete set
	\[
	U = \{ x \in X \mid \{x\} \text{ is closed and open}\}.
	\]
	Assume that $X\setminus U\neq\emptyset$. By scatteredness, $\{x\}$ is open in $X \setminus U$ for some $x \in X \setminus U$.
	Therefore $X \setminus (U \cup \{x\})=(X \setminus U) \setminus \{x\}$ is closed in $X \setminus U$ and hence closed in $X$.
	Thus both $\{x\}$ and $X \setminus (U \cup \{x\})$ are closed subsets. Since $X$ is compact Hausdorff, there are disjoint open subsets $V_1$ and $V_2$ containing $x$ and $X \setminus (U \cup \{x\})$, respectively. We may assume $V_1$ is closed because $X$ is Stonean.
	Observe that $V_1$ is infinite; otherwise $V_1 \setminus \{x\}$ would be closed and $\{x\}=V_1 \setminus (V_1 \setminus \{x\})$ open, contradicting $x \notin U$.
	Hence $V_1 \setminus \{x\}$ is infinite, too.
	Pick two disjoint infinite subsets $W_1,W_2$ covering $V_1 \setminus \{x\}$. Since $W_i$ is contained $U$, it must be open.
	If $W_i$ were closed, then it is compact, contradicting that it is both discrete (as subset of $U$) and infinite. So $W_i \subsetneq \overline{W_i} \subseteq V_1$.
	Moreover, $\overline{W_i} \cap W_j = \emptyset$ for $i \neq j$, and $W_1 \cup W_2 = V_1 \setminus \{x\}$, so $\overline{W_i} = W_i \cup \{x\}$ whence $\overline{W_1} \cap \overline{W_2} = \{x\}$. Since $X$ is Stonean, $\overline{W_i}$ is open, whence $\{x\}$ is open, contradicting $x\notin U$. Hence $X=U$, and since $U$ is discrete and compact, it must be finite.
\end{proof}

\begin{theorem}\label{thm:AW}
	The following are equivalent for an AW*-algebra $A$:
	\begin{itemize}
		\item[(1)] $A$ is finite-dimensional;
		\item[(2)] $\CC(A)$ is algebraic;
		\item[(3)] $\CC(A)$ is continuous;
		\item[(4)] $\A(A)$ is algebraic;
		\item[(5)] $\A(A)$ is continuous.
	\end{itemize}
\end{theorem}
\begin{proof}
	Clearly, if $A$ is finite-dimensional, it is scattered, which we know to be equivalent with $\CC(A)$ being algebraic (cf. Theorem~\ref{thm:algebraic}) or continuous (cf. Theorem~\ref{thm:continuous}). Conversely, assume that $A$ is scattered. Theorem~\ref{thm:Kusuda} then implies that all maximal commutative C*-subalgebras of $A$ are scattered. But maximal C*-subalgebras are automatically AW*-algebras by Proposition~\ref{prop:galois},
	and scattered commutative AW*-algebras are finite-dimensional by Lemma~\ref{lem:stonescattered}.
	Since all maximal commutative C*-subalgebras of $A$ are finite-dimensional, so is $A$ itself~\cite[Exercise 4.6.12]{kadisonringrose:oa1}.
	It follows that (1), (2), and (3) are equivalent. Moreover, since the class of finite-dimensional C*-algebras coincides with the class of finite-dimensional AW*-algebras, it follows that $\A(A)=\CC(A)$ if $A$ is finite-dimensional. Hence if $A$ is finite-dimensional, then $\A(A)$ is algebraic, hence also continuous since all algebraic dcpos are continuous \cite[Proposition I-4.3]{gierzetal:domains}. So (1) implies (4), which implies (5). 
	
	Finally, we show that (5) implies (1) by contraposition. Suppose $A$ is infinite-dimensional. Pick a maximal commutative C*-subalgebra $C \subseteq A$; its Stonean spectrum $X$ will have infinitely many points, and by compactness a non-isolated point $x$.
	Any other point $y_1 \in X$ is separated from $x$ by a clopen $U_1$, and induction gives a sequence of disjoint clopens $U_1,U_2,\ldots$. Their indicator functions form an infinite set $P$ of pairwise orthogonal projections in $A$.
	
	Let $I \subseteq P$ be an infinite subset with infinite complement. Its supremum $p=\bigvee I$ is nonzero. Choosing some nonzero $q \in P \setminus I$ gives $rq=0$ for each $r \in I$, and hence $pq=0$~\cite[Proposition~3.6]{berberian}, so that $p \neq 1$. By Lemma~\ref{lem:atomsinCA}, $C^*(p)$ is an atom in $\CC(A)$.
	Consider the directed family $\{ C^*(F) \mid F \subseteq P \text{ finite} \}$ of elements of $\A(A)$, whose supremum contains $p$. We will show that no element $C=C^*(p_1,\ldots,p_n)$ of the family can contain $p$. Observe that $p_{n+1}=1-\sum_{i=1}^n p_i \in C$ is orthogonal to each $p_1,\ldots,p_n$, and hence $\sum_{i=1}^{n+1}p_i=1$. Therefore $C = C^*(p_1,\ldots,p_{n+1}) = \Span\{p_1,\ldots,p_{n+1}\} = \Span\{p_1,\ldots,p_n,1\}$.
	If $p$ were in $C$, we could thus write it as $p=\sum_{i=1}^n \lambda_i p_i + \lambda 1$ for some coefficients $\lambda,\lambda_i$.
	Pick a nonzero element $q \in I$ distinct from $p_1,\ldots,p_n$. Because $q \leq p$ we find $q = qp = \sum_{i=1}^n\lambda_iqp_i + \lambda q = \lambda q$, whence $\lambda=1$.
	Now pick a nonzero element $q \in P \setminus I$ distinct from $p_1,\ldots,p_n$. Then $qp=0$ and hence $qp_i=0$ for each $i=1,\ldots,n$.
	Thus $q = \sum_{i=1}^n \lambda_i qp_i + q = qp = 0$, and $p$ cannot be contained in $C$.
	Therefore $C^*(p)$ is not compact, but since it is an atom of $\A(A)$, it is way above the bottom element only. 
	Hence $\bigvee \ddown C^*(p) \neq C^*(p)$, and $\A(A)$ is not continuous.
\end{proof}

We conclude that, at least from a domain-theoretic perspective on $\CC(A)$, C*-algebras are more interesting than W*-algebras or AW*-algebras, since C*-algebras contain a subclass of infinite-dimensional algebras $A$ for $\CC(A)$ is a domain, whereas the only subclass of AW*-algebras or of W*-algebras for which $\CC(A)$ is a domain is the class of finite-dimensional algebras. Nevertheless, the most satisfactory results for reconstructing the structure of $A$ from $\CC(A)$ are obtained for AW*-algebras and W*-algebras, as discussed in the Introduction.
Any C*-algebra $A$ can be extended to a W*-algebra by taking its double dual $A^{**}$, also called the \emph{enveloping W*-algebra}~\cite[\S III.2]{takesaki:oa1}. This in fact gives an adjunction of categories showing that W*-algebras form a non-full reflective subcategory of C*-algebras~\cite[3.2]{dauns:tensorproduct}.
There are many examples of C*-algebras $A$ for which $\CC(A)$ is continuous but $\CC(A^{**})$ is not: any infinite-dimensional scattered C*-algebra will do, such as $C(X)$ for the infinite compact Hausdorff scattered spaces $X$ of Example~\ref{ex:scattered}.

\section{Directed colimits}\label{sec:colimits}

 Because an AF-algebra is a directed colimit of finite-dimensional C*-algebras, one might wonder whether the functors $\CC$ or $\CC_\AF$ preserve directed colimits of C*-algebras. In general the answer for both functors is negative. It turns out to be useful to treat the case of $\CC_\AF$ first, for which we need to calculate the directed colimit of algebraic dcpos. The answer for $\CC$ can then be derived. Since $\CC_\AF\simeq\B\circ\Proj$, it is useful to split the case of $\CC_\AF$ in two other cases, namely the behaviour $\Proj$ with respect of directed colimits of C*-algebras, and the behaviour of $\B$ with respect to directed colimits of orthomodular posets. It turns out that $\B$ does preserve directed colimits, but $\Proj$ does not, from which it follows that $\CC_\AF$ does not preserve directed colimits either.
  
Let $A=\overline{\bigcup_{i\in I}A_i}$ where $\{A_i\}_{i\in I}$ is a directed collection of C*-subalgebras. 
When $i \leq j$, the inclusion $A_i\subseteq A_j$ makes $\Proj(A_i)\subseteq\Proj(A_j)$ into a sub-orthomodular poset. More generally, consider a collection $\{P_i\}_{i\in I}$, where $I$ is directed, and $P_i$ is a sub-orthomodular poset of $P_j$ if $i\leq j$. The colimit $P$ of the $P_i$ in the category of orthomodular posets and orthomodular morphisms is~\cite[Theorem 4.10]{navararogalewicz:pasting}
\[P=\bigcup_{i\in I}P_i\text,\]
where $P$ is ordered by $p\leq q$ if and only if there is some $i\in I$ such that $p\leq q$ in $P_i$. The orthocomplementation of $P$ is given by $p^\perp=q$ if and only if there is some $i\in I$ such that $p^\perp=q$ in $P_i$. 
  
Proposition \ref{prop:Bisfunctor} shows the inclusions $P_i\subseteq P_j$ induce Scott continuous inclusions $\B(P_i)\subseteq\B(P_j)$, and Proposition \ref{prop:BofPisalgebraiccompletesemilattice} shows that the compact elements of $\B(P_i)$ and $\B(P_j)$ are their finite Boolean subalgebras, so $\K(\B(P_i))\subseteq\K(\B(P_j))$. To compute $\colim_{i\in I}\B(P_i)$, we need a description of the directed colimit of algebraic dcpos $\{X_i\}_{i\in I}$ such that $\K(X_i)\subseteq\K(X_j)$ if $i\leq j$. This description seems to be folklore; for the sake of completeness we include a proof. We first recall a definition and two lemmas about free directed completions.  
  
  \begin{definition}
  	Let $X$ be a poset. Then we denote the set of all \emph{ideals} of $X$, i.e., all directed subsets $I\subseteq X$ such that $I=\down I$, by $\Id(X)$. 
  \end{definition}

  \begin{lemma}\cite[Proposition I-4.10]{gierzetal:domains}\label{lem:idealcompletion}
  	Let $X$ be a poset. Then $\Id(X)$ is an algebraic dcpo ordered by inclusion, where unions provide directed suprema. Moreover, $x \mapsto \down x$ is an order embedding $X\mapsto\Id(X)$, $x\mapsto\down x$ with image $\K(\Id(X))$, the set of all compact elements of $\Id(X)$. If $X$ itself is an algebraic dcpo, then $X\simeq\Id(\K(X))$.	
  \end{lemma}
  
  \begin{lemma}\cite[Corollary 3.1.6]{stoltenberghansen:domains}\label{lem:extension}
  	Let $X$ and $Y$ be dcpos with $X$ algebraic. Any monotone map $\varphi \colon \K(X)\to Y$ has a unique Scott continuous extension $X\to Y$ given by \[x \mapsto \bigvee\{\varphi(c):c\in\K(X)\cap\down x\}. \]
  \end{lemma}

  \begin{proposition}\label{prop:directedcolimofalgebraicdcpos}
  	If $\{X_i\}_{i\in I}$ be a directed family of algebraic dcpos such that $X_i\subseteq X_j$ and $\K(X_i)\subseteq\K(X_j)$ when $i \leq j$, then
  	\[\colim_{i\in I}X_i=\Id\left(\bigcup_{i\in I}\K(X_i)\right)\text.\] 
  	The colimiting cone $\varphi_i \colon X_i \to\colim_{i\in I}X_i$ is given by
    $
      x \mapsto\bigcup\{\down c \mid c\in\K(X_i)\cap\down x\}$.
  \end{proposition}	
  \begin{proof}
  	Write $X$ for the proposed colimit. Lemma \ref{lem:idealcompletion} makes it an algebraic dcpo whose compact elements are precisely the elements $\down c$, where $c$ is compact in $X_i$ for some $i\in I$. We have to show that the $\varphi_i$ are well-defined Scott continuous functions such that the upper triangle of the following diagram commutes:
  	
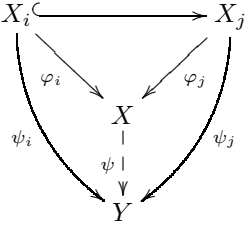
\begin{align*}\xymatrix{X_i\ar@{^{(}->}[rr]\ar[dr]_{\varphi_i}\ar@/_1pc/[ddr]_{\psi_i} && X_j\ar[dl]^{\varphi_j}\ar@/^1pc/[ddl]^{\psi_j}\\
  		& X\ar@{-->}[d]_\psi &\\
  		& Y & }
  	\end{align*}
  	Let $\theta_i \colon \K(X_i)\to X$ be the map $c\mapsto\down c$. Then $\theta_i$ is monotone, since it is the restriction of the canonical embedding in Lemma \ref{lem:idealcompletion}. Hence it has a unique Scott continuous extension by Lemma \ref{lem:extension}, which is clearly equal to $\varphi_i$. Therefore $\varphi_i$ is a well-defined Scott continuous map. To see that $\{\varphi_i\}$ forms a cocone, we have to show that $\varphi_i(x)=\varphi_j(x)$ for each $x\in X_i$ if $i\leq j$. Since $\K(X_i)\subseteq\K(X_j)$ by assumption, it follows that $\theta_j$ is an extension of $\theta_i$. Thus the restriction of $\varphi_j$ to $\K(X_i)$ equals $\theta_i$, and Lemma \ref{lem:extension} guarantees that $\varphi_i$ equals the restriction of $\varphi_j$ to $X_i$. 
  	
  	Now let $Y$ be another dcpo and $\psi_i \colon X_i\to Y$ Scott continuous maps such that $\psi_i(x)=\psi_j(x)$ for each $x\in X_i$ and $i\leq j$. We prove that there is a unique Scott continuous map $\psi \colon X\to Y$ making the diagram above commute. Let $\eta \colon \bigcup_{i\in I}\K(X_i)\to Y$ be the map $c\mapsto\psi_i(c)$ if $c\in X_i$. Since $\K(X_i)\subseteq\K(X_j)$ and $\psi_j|_{X_i}=\psi_i$ if $i\leq j$, this map is a well-defined monotone map. By Lemma \ref{lem:idealcompletion}, the map $\kappa \colon \bigcup_{i\in I}\K(X_i)\to\K(X)$ given by $c\mapsto\down c$ is an order isomorphism, so $\eta\circ\kappa^{-1} \colon \K(X)\to Y$ is monotone, and so Lemma \ref{lem:extension} guarantees a unique Scott continuous extension $\psi \colon X\to Y$. 
  	Let $c\in \K(X_i)$. Then $c\in\bigcup_{i\in I}\K(X_i)$, so $\kappa(c)=\down c$, whence
  	\[\psi\circ\varphi_i(c)=\psi\circ\theta_i(c)=\psi(\down c)=\eta\circ\kappa^{-1}(\down c)=\eta(c)=\psi_i(c)\text.\]
  	So the restriction of $\psi\circ\varphi_i$ to $\K(X_i)$ equals the restriction of $\psi_i$ to $\K(X_i)$. Lemma \ref{lem:extension} now gives $\psi\circ\varphi_i=\psi_i$. Now assume that $\psi' \colon X\to Y$ is another Scott continuous map making the diagram commute, and let $J\in\K(X)$. By definition, $J=\kappa(c)$ for some $i\in I$ and $c\in\K(X_i)$. Since $c\in\K(X_i)$, 
  	\[\varphi_i(c)=\bigcup\{\down c' \mid c'\in\K(X_i)\cap\down c\}=\down c=\kappa(c)\text,\] whence
  	\[\psi'(J)=\psi'\circ\kappa(c)=\psi'\circ\varphi_i(c)=\psi_i(c)=\eta(c)=\eta\circ\kappa^{-1}(J)=\psi(J),\]
  	so $\psi'$ coincides with $\psi$ on $\K(X)$. Now $\psi'=\psi$ by Lemma \ref{lem:extension}.
  \end{proof}

  \begin{corollary}
  	The functor $\CC$ does not preserve directed colimits.
  	\end{corollary}
  \begin{proof}
    Let $X$ be the Cantor space and $A=C(X)$. This is a commutative AF-algebra that is not scattered. So $A$ is the directed colimit of a set $\{A_i\}_{i\in I}$ of finite-dimensional C*-subalgebras, i.e.\ $A=\overline{\bigcup_{i\in I}A_i}$. 
    If $i\leq j$, there is an inclusion $f \colon A_i\to A_j$. 
    Now $\CC(A_i) \subseteq \CC(A_j)$ because $\CC(f)$ is an order embedding (by the remark below Proposition \ref{prop:CAisdcpo}). By Proposition \ref{prop:finitedimiscompact} then $\K(\CC(A_i))=\CC(A_i)\subseteq\CC(A_j)=\K(\CC(A_j))$, and the previous proposition shows that the colimit of the $\CC(A_i)$ exists and is algebraic. However, since $A$ is not scattered, $\CC(A)$ is not algebraic, and so $\CC(A)\not\simeq\colim_{i\in I}\CC(A_i)$.
 	\end{proof}

Next we show that the functors $\A$ and $\V$ do not preserve directed colimits either. As morphisms of AW*-algebras we choose $*$-homomorphisms that are \emph{normal}, i.e. that preserve suprema of projections. W*-algebras form a full subcategory, so finite-dimensional C*-algebras are certainly AW*-algebras. We first have to show the existence of the directed colimit of some diagram of finite-dimensional C*-algebras in the category of AW*-algebras. To do so, let $A_n=\C^{n+1}$, and define for each $n\in\omega$ the function $f_{n,n+1}\colon A_n\to A_{n+1}$ by $(x_0,x_1,\ldots,x_n)\mapsto(x_0,x_1,\ldots,x_n,x_n)$, and for $m<n$, define $f_{m,n}\colon A_m\to A_n$ by \[f_{m,n}=f_{n-1,n}\circ f_{n-2,n-1}\circ\cdots \circ f_{m+1,m+2}\circ f_{m,m+1}.\] 
It is easy to see that $f_{m,n}$ is a $*$-homomorphism. Note that the choice of $A_n=\C^{n+1}$ instead of $A_n=\C^n$ assures that also $f_{0,1}$ is unital, as required in our definition of a $*$-homomorphism. Since the domain $\C^{m+1}$ of $f_{m,n}$ is finite-dimensional, and commutative, it only has a finite number of projections, namely the elements $(a_0,a_1,\ldots,a_n)$ with $a_i\in\{0,1\}$, and so $f_{m,n}$ must be normal. It follows that the $f_{m,n}$ form a directed diagram in both the category of AW*-algebras and in its full subcategory of W*-algebras. Recall the definition of $\ell^\infty(X)$ from Example \ref{ex:l infty}. 
The next lemma shows that $\ell^\infty(\omega+1)$ is the colimit of the directed diagram defined above. Since $\omega$ and $\omega+1$ have the same cardinality, we could have considered $\ell^\infty(\omega)$ instead, but $\ell^\infty(\omega+1)$ is handier for notational reasons. For each $n\in\omega$, define $g_n \colon A_n\to \ell^\infty(\omega+1)$ by 
\[ 
  g_n(x_0,x_1,\ldots,x_n)(i)=\begin{cases} 
    x_i & i\leq n;\\
    x_n &i\geq n, 
  \end{cases}
\]
which is a normal $*$-homomorphism for the same reasons as $f_{m,n}$.
Since clearly $g_{n+1}\circ f_{n,n+1}= g_n$, so that $g_n\circ f_{m,n}=g_m$ for each $n\geq m$ in $\omega$, and the $g_n$ form a cocone. 

\begin{lemma}\label{lem:reviewer1}
  The cocone $(\ell^\infty(\omega+1),g_n)$ is a colimit of diagram $\{f_{m,n} \colon A_m \to A_n\}_{m,n\in\omega}$, both in the category of AW*-algebras with normal $*$-homomorphisms, and in its full subcategory of W*-algebras.
\end{lemma}
\begin{proof}
	Let $A$ be an AW*-algebra, and for each $n\in\omega$, let $h_n\colon A_n\to A$ be a normal $*$-homomorphism satisfying $h_n\circ f_{m,n}=h_m$ for each $m<n$ in $\omega$.  
	We have to establish a unique normal $*$-homomorphism $k\colon \ell^\infty(\omega+1)\to A$ such that $k \circ g_m = h_m$ for all $m$.
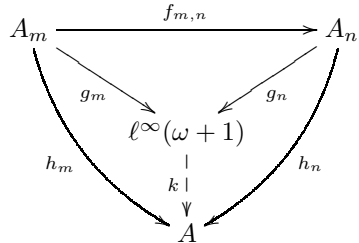
		\begin{align*}\xymatrix{A_m\ar[rr]^{f_{m,n}}\ar[dr]_{g_m}\ar@/_1pc/[ddr]_{h_m} && A_{n}\ar[dl]^{g_n}\ar@/^1pc/[ddl]^{h_n}\\
			& \ell^\infty(\omega+1)\ar@{-->}[d]_k &\\
			& A & }
		\end{align*}
	For each $0\leq i\leq \omega$, let $e_i\colon\omega+1\to\C$ be given by $e_i(j)=\delta_{ij}$. Note that any projection in $\ell^\infty(\omega+1)$ is the characteristic function $\chi_S$ of some subset $S$ of $\omega+1$, so $S\mapsto\chi_S$ is a Boolean isomorphism of the power set of $\omega+1$ to $\Proj(\ell^\infty(\omega+1))$. Since power sets are atomistic Boolean algebras, where the singletons are the atoms, it follows that $e_i=\chi_{\{i\}}$ is an atomic projection, and $\chi_S=\bigvee_{i\in S}e_i$.

	Similarly, let $e_i^n\in A_n$ be the projection given by $(e_i^n)_j=\delta_{ij}$. Notice $g_{n+1}(e^{n+1}_n)=e_n$ for each $n\in\omega$. Next, for each $n\in\omega$ set $p_n=h_{n+1}(e_n^{n+1})$. Now $p_n$ is a projection in $A$ because $h_{n+1}$ is a $*$-homomorphism. If $m\neq n$, say $m<n$, then
	\begin{align*}
	  p_np_m & =h_{n+1}(e_n^{n+1})h_{m+1}(e_m^{m+1})=h_{n+1}(e_n^{n+1}) \cdot h_{n+1}(f_{m+1,n+1}(e_m^{m+1}))\\
	  & =h_{n+1}(e_n^{n+1})h_{n+1}(e_m^{n+1})=h_{n+1}(e^{n+1}_ne^{n+1}_m)=0,  
	  \end{align*}
	since $e^{n+1}_n$ and $e^{n+1}_m$ are orthogonal. 
	Let $p_\omega=\left(\bigvee_{n\in\omega}p_n\right)^\perp$. Clearly $p_\omega$ is orthogonal to any other $p_n$, and moreover, $\bigvee_{n\in\omega+1}p_n=1$.

	Now, define $k$ as follows. Let $e\in \ell^\infty(\omega+1)$ be a projection, so $e=\chi_S$ for some $S\subseteq \omega+1$. Then we define $k(e)=\bigvee_{i\in S}p_i$. Since $\bigvee_{i\in\omega+1} p_i=1$, it follows that $k(1)=1$, which is the reason to consider $\omega+1$ instead of $\omega$. Let $e_1$ and $e_2$ be projections in $\ell^\infty(\omega+1)$, say $e_1=\chi_S$ and $e_2=\chi_T$ for some $S,T\subseteq\omega+1$. Then $e_1e_2=\chi_S\chi_T=\chi_{S\cap T}$, so $\Proj(\ell^\infty(\omega+1))$ is closed under multiplication. Moreover, since the $p_i$ are mutually orthogonal,
	\[ g(e_1)g(e_2)=\bigvee_{i\in S}p_i\bigvee_{j\in T}p_j=\bigvee_{i\in S}\bigvee_{j\in T}p_ip_j=\bigvee_{i\in S\cap T}p_i=g(\chi_{S\cap T})=g(e_1e_2).\]
	Next extend $k$ linearly to the span of $\Proj(\ell^\infty(\omega+1))$. The extended function then is linear and still multiplicative. Moreover, since projections are self-adjoint, $k$ preserves the involution. Since $\ell^\infty(\omega+1)$ is a commutative AW*-algebra, it is generated by its projections. In a commutative C*-algebra the product of two projections is again a projection, so the span of the projections in $\ell^\infty(\omega+1)$ is dense.
	Thus $k$ extends uniquely to a $*$-homomorphism $\ell^\infty(\omega+1)\to A$, which is furthermore normal because its restriction to projections by definition preserves suprema.
	Now, for each $n\in\omega$ and $1\leq i < n$,
	\[  k\circ g_n(e^{n}_i)=k(e_i)=p_i=h_{i+1}(e^{i+1}_i) = h_{n}\circ f_{i+1,n}(e^{i+1}_i) = h_n(e^n_i).  \]
	Since also $k\circ g_n(1)=1=h_n(1)$, and $A_n$ is spanned by $e_0^n,e_1^n,\ldots,e_{n-1}^n,1$, it follows that $k\circ g_n=h_n$, establishing existence of $k$.

	If $k'\colon \ell^\infty(\omega+1)\to A$ is another normal $*$-homomorphism satisfying $k'\circ g_n=h_n$, then
	\[ p_n=h_{n+1}(e^{n+1}_n)=k'\circ g_n(e^{n+1}_n)=h_n(e_n).\]
	Moreover, since $k'$ is normal, for any projection $e=\chi_S$ in $l^{\infty}(\omega+1)$,
	\[ k'(e)=k'(\chi_S)=k'\left(\bigvee_{i\in S}e_i\right)=\bigvee_{i\in S}k'(e_i)=\bigvee_{i\in S}p_i=k(e),\]
 	making $k'$ coincides with $k$ on $\Proj(\ell^\infty(\omega+1))$. As we uniquely extended $k$ to a linear map on the span of $\Proj(\ell^\infty(\omega+1))$, it follows that $k$ and $k'$ coincide on the span of the projections of $\ell^\infty(\omega+1)$. Since $k$ and $k'$ are both continuous, they also coincide on $\ell^\infty(\omega+1)$. 
\end{proof}

\begin{proposition}\label{prop:AandVdontpreservecolims}
  The functors $\A$ and $\V$ do not preserve directed colimits.
  \end{proposition}
  \begin{proof}
  By Lemma \ref{lem:reviewer1}, we can consider a directed set $\{A_i\}_{i\in I}$ of finite-dimensional C*-algebras for which the colimit exists both in the category of AW*-algebras with normal $*$-homomorphisms, and in the full subcategory of W*-algebras, and such that the colimit is infinite-dimensional. Since $\A(A_i)=\V(A_i)$ is algebraic by Theorem \ref{thm:AW}, it follows from Proposition \ref{prop:directedcolimofalgebraicdcpos} that $\colim_{i\in I}\A(A_i)$ is algebraic, too. However, since $A$ is infinite-dimensional, it follows from Theorem \ref{thm:AW} that $\A(A)=\V(A)$ cannot be algebraic. Hence $\A$ and $\V$ cannot preserve colimits.
\end{proof}

In contrast to $\CC$, $\A$, and $\V$, the functor $\B$ does preserve directed colimits.

\begin{proposition}\label{prop:Bpreservesdirectedcolimits}
  The functor $\B \colon \OMP\to\DCPO$ preserves directed colimits. 
\end{proposition}
\begin{proof}
  Let $(\{P_i\}_{i\in I},\{\varphi_{ij} \colon P_i\to P_j\}_{i,j\in I})$ be a directed system in $\OMP$ with colimit $P$. Write $\varphi_i \colon P_i \to P$ for the colimiting cone. Replacing $P_i$ by $\varphi_j[P_i]$ if necessary, we may assume that the $\varphi_{ij}$ are inclusions and $P=\bigcup_{i\in I}P_i$. For each orthomodular poset $Q$, the compact elements of $\B(Q)$ are precisely the finite Boolean subalgebras of $Q$ by Proposition \ref{prop:BofPisalgebraiccompletesemilattice}. Hence $\K(\B(P_i))\subseteq\K(B(P_j))$ if $i\leq j$, and Proposition \ref{prop:directedcolimofalgebraicdcpos} shows $\colim_i \B(P_i)=\Id(\bigcup_i\K(\B(P_i))$. Since $P_i$ embeds into $P$, it follows that every finite Boolean subalgebra of $P_i$ is a finite Boolean subalgebra of $P$, whence $\bigcup_i\K(\B(P_i))\subseteq\K(P)$. Conversely, if $B=\{b_1,\ldots,b_n\}\subseteq P$ is a finite Boolean subalgebra, there is ${i_k}\in I$ such that $b_k\in P_{i_k}$ for $k\in\{1,\ldots,n\}$. Directedness now produces $i \geq i_1,\ldots,i_n$, and $B\subseteq P_i$. 
  	
 	Next we show that $B$ is a Boolean subalgebra of $P_i$ by showing that it is a sub-orthomodular poset. Let $b\in B$, and let $b^\perp$ be its orthocomplement in $B$. Since $B$ is a Boolean subalgebra, $b^\perp$ is also the orthocomplement of $b$ in $P$. Since $P_i$ is a sub-orthomodular poset, $b^\perp$ is also the orthocomplement of $b$ in $P_i$, and $B$ is closed under the orthocomplementation of $P_i$. Let $a,b\in B$ be orthogonal in $B$. Then $a\leq b^\perp$ in $B$ and therefore in $P$. The join $a\vee_Pb$ of $a$ and $b$ in $P$ lies in $B$, hence must equal the join $a\vee_Bb$ of $a$ and $b$ in $B$. Since $a\leq b^\perp$ in $P$ and $a,b\in P_i$, also $a\leq b^\perp$ in $P_i$. It follows that $a\vee_Pb\in P_i$, whence the join $a\vee_{P_i}b$ of $a$ and $b$ in $P_i$ exists and equals $a\vee_{P_i}b$. Consequently $a\vee_{P_i}b=a\vee_Bb$, and $B$ is a sub-orthomodular poset of $P_i$. Since it is a finite Boolean algebra, $B\in\K(\B(P_i))$.
  	
 	We conclude that $\K(\B(P))\subseteq\bigcup_{i\in I}\K(\B(P_i))$, hence $\K(\B(P))=\bigcup_{i\in I}\B(P_i)$. It now follows from Lemma \ref{lem:idealcompletion} that
  	\[\B(P)\simeq\Id(\K(\B(P)))\simeq\Id\left(\bigcup_{i\in I}\K(\B(P_i))\right)=\colim_{i\in I}\B(P_i),\] which is exactly what we wanted to prove.  	  
  \end{proof}	
  
We return to the case where $P_i=\Proj(A_i)$ for some directed set $\{A_i\}_{i\in I}$ of C*-subalgebras of a C*-algebra $A$ that contains $\bigcup_{i\in I}A_i$ as a subset, so that $A=\colim_{i\in I}A_i$. Say that $A$ has the \emph{lattice property} if $\Proj(A)$ is a lattice, and that it has the \emph{directed set property} if its collection of finite-dimensional C*-algebras is directed. 
In case $A$ is approximately finite-dimensional, the lattice property and the directed set property can be related to each other.

\begin{lemma}\label{lem:dspandlp}
	Let $A=\overline{\bigcup_{i\in I}A_i}$ be an AF-algebra, where $\{A_i\}_{i\in I}$ is some directed collection of finite-dimensional C*-subalgebras of $A$. 
	\begin{enumerate}
		\item[(a)] If $\Proj(A)\simeq\bigcup_{i\in I}\Proj(A_i)$, then $\Proj(A)$ is a lattice;
		\item[(b)] If $A$ has the directed set property and $\{A_i\}_{i\in I}$ is the collection of \emph{all} finite-dimensional C*-algebras, then $\Proj(A)\simeq\bigcup_{i\in I}\Proj(A_i)$. 
	\end{enumerate}
\end{lemma}
\begin{proof}
	Let $p,q\in\bigcup_{i\in I}\Proj(A_i)$. Then $p\in A_i$ and $q\in A_j$ for some $i,j\in I$. Directedness gives $k\geq i,j$. Since $A_k$ is finite-dimensional, $r = p \vee q$ exists in $A_k$. Let $s\in\bigcup_{i\in I}\Proj(A_i)$ be an upper bound of $p$ and $q$, say $s\in A_m$. Again directedness gives $n\geq k,m$. Now, $A_k\subseteq A_n$ are both finite-dimensional and so W*-algebras, and $A_k$ is a W*-subalgebra of $A_n$. Therefore the join of $p$ and $q$ in $A_n$ equals their join in $A_k$. It follows that $r$ is also the join of $p$ and $q$ in $A_n$, and $r\leq s$. Thus $r = p \vee q$ in $A$. Similarly, $p \wedge q$ exists in $\bigcup_{i\in I}\Proj(A_i)$ and hence in $\bigcup_{i\in I}\Proj(A_i)$, and $\Proj(A)$ is a lattice.
	
	For (b), we have $A_i\subseteq A$ hence $\Proj(A_i)\subseteq\Proj(A)$ for each $i\in I$, so $\bigcup_{i\in I}\Proj(A_i)\subseteq\Proj(A)$. Let $p\in\Proj(A)$. Then $p\in C^*(p)$, which is finite-dimensional by Lemma \ref{lem:atomsinCA}, hence $p\in A_j$ for some $j\in I$. It follows that $p\in\Proj(A_j)$, whence $p\in\bigcup_{i\in I}\Proj(A_i)$.
\end{proof}

Combining both statements in the previous lemma shows that an AF-algebra with the directed set property has the lattice property. This has been shown by Lazar in \cite[Theorem 3.4]{lazar:afalgebraslattice}, who also showed the converse for separable AF-algebras. It is remarkable that the fact that there is \emph{some} directed set of finite-dimensional C*-subalgebras whose union is dense in $A$ does not imply that \emph{all} finite-dimensional C*-subalgebras of $A$ are directed. This follows from Lazar's construction of a separable AF-algebra $A$ without the lattice property \cite{lazar:afalgebrasdirected}.

\begin{corollary}
	The functors $\Proj \colon \CStar\to\OMP$ and $\CC_{\AF} \colon \CStar \to \DCPO$ do not preserve directed colimits.
\end{corollary}
\begin{proof}
  Let $A$ be Lazar's AF-algebra, and let $\{A_i\}_{i\in I}$ be the directed set of finite-dimensional C*-subalgebras of $A$ such that $A=\overline{\bigcup_{i\in I}A_i}$.  
  Since $A$ does not have the lattice property, it follows from Lemma \ref{lem:dspandlp} that $\Proj(A)\not\simeq\colim_{i\in I}\Proj(A_i)$, where we used that the directed colimit of orthomodular posets is given by their union. We conclude that $\Proj$ does not preserve directed colimits. Then we also have
  \[\colim_{i\in I}\B\circ\Proj(A_i)\simeq \B\left(\colim_{i\in I}\Proj(A_i)\right)\not\simeq\B\circ\Proj\left(\colim_{i\in I}A_i\right),\]
  where the equivalence follows from Proposition \ref{prop:Bpreservesdirectedcolimits}, and the inequivalence follows from the statement that $\Proj$ does not preserve directed colimits combined with the fact that the functor $\B$ determines orthomodular posets up to isomorphism as proven in \cite{hardingheunenlindenhoviusnavara:booleansubalgebras}. We conclude that $\B\circ\Proj$ cannot preserve directed colimits. Since $\CC_{\AF} \simeq \B \circ \Proj$ by Theorem \ref{thm:CAFisBProjA}, it follows that neither $\CC_{\AF}$ can preserve directed colimits.
\end{proof}

\section*{Acknowledgement}

This article extends the earlier conference proceedings version~\cite{heunenlindenhovius:lics}. Chris Heunen was supported by the Engineering and Physical Sciences Research Council Fellowship EP/L002388/1. 
Bert Lindenhovius was supported by the U.S.\ Department of Defence and the U.S. Air Force Office of Scientific Research under the MURI grant number FA9550-16-1-0082
entitled, ``Semantics, Formal Reasoning, and Tool Support for Quantum
Programming'', and by the Netherlands Organisation for Scientific Research under TOP-GO grant 613.001.013.
We are indebted to Klaas Landsman for useful comments, John Harding for pointing out Lemma~\ref{lem:equivalencerelations}, Michael Mislove for helping us prove Lemma~\ref{lem:orderscattered}, and for pointing out how to compute a colimit of algebraic dcpos, and Ruben Stienstra for helping us prove Lemma~\ref{lem:stonescattered}.
Finally, we are very grateful to anonymous reviewers for many detailed and insightful comments, and in particular for pointing out and proving Lemma~\ref{lem:reviewer1} and Proposition~\ref{prop:AandVdontpreservecolims}.

\bibliographystyle{apalike}
\bibliography{domains}

\end{document}